\pdfoutput=1
\documentclass[11pt]{article}
\usepackage{amssymb}
\usepackage{graphicx}

\parskip \baselineskip

\newcommand{\ra}{\Rightarrow}
\newcommand{\s}{\subseteq}
\newcommand{\es}{\emptyset}

\newcommand{\La}{\mathbb{L}}

\newcommand{\Pow}{\mathbb{P}}

\newcommand{\Z}{\mathbb{Z}}
\newcommand{\B}{{\cal B}}
\newcommand{\ol}{\overline}
\newcommand{\ul}{\underline}
\newcommand{\Ra}{\Rightarrow}
\newcommand{\LRa}{\Leftrightarrow}

\newcommand{\LL}{\Lambda}
\newcommand{\vv}{\vert}

\begin{document}
\centerline{\bf Modular lattices of finite length (Part B)}
\bigskip
\centerline{Marcel Wild}\bigskip

By definition 'Part B' means Chapters 7 to 9 of a larger project that is outlined in "Modular lattices of finite length (Part A)".
The overall project is planned to consist of these Sections:

\begin{enumerate}
	\item Introduction
	\item Basic facts
	\item Lattice congruences via  coverings
	\item Geometric lattices and matroids
	\item Partial linear spaces and matroids
	\item Enumeration of modular lattices\\
	\item Bases of lines: The essentials
	\item Calculating the submodule lattice of a finite $R$-module
	\item Bases of lines: Localization and two types of cycles\\
	\item Existence of R-linear representations, partition representations, and 2-distributivity
	\item Classification of the k-linear representations of acyclic modular lattices
	\item Axiomatization\\
	\item Constructions: Gluings and subdirect products
	\item Relatively free structures: From semilattices to distributive to modular lattices
	
\end{enumerate}

\vskip 3cm
\centerline{\bf 7. Bases of lines: The essentials}
\medskip

Inspired by projective spaces (Sec. 5.6.2) let us jump into medias res and define:

 \begin{enumerate}
 	\item [(1)] A \underbar{line} of any\footnote{We will soon return to the modular case but for the time being let us see how far we get without that restriction.}  lattice $L$ is
subset $l\s J(L)$ with $\vv l\vv\ge 3$ and
maximal with the {\it same-join-property}, i.e. any $p\not =q$ in $l$ yield the
same join, which we denote by $\overline{l}$.
\end{enumerate}

Consider $L=SM_{10}$ in Figure 2.1. It has the line $\ell=\{2,6,8\}$ because\footnote{Do not show these equalities to your 8-year old youngster.} $2+6=2+8=6+8=1$, and neither the join-irreducibles 3 nor 4 can be added without destroying the same-join-property.

The following gadget can fabricate lines.  A length two interval $[x_0,x]\simeq M_n$
of a lattice $L$ is a \underbar{line-interval} if 
$[x_0,x]$ contains {\it all} lower covers of $x$ in $L$. Recall that $n\ge 3$ by definition of $M_n$. Thus $[x_0,x]=[10,18]$ in the lattice $L_1$ of Figure 3.1 is a line-interval, but not $[x_0,x]=[3,1]$ in $SM_{10}$ because it does not contain the lower cover $7\prec x$.

\smallskip
{\bf Lemma 7.1 [HW,2.3]: }{\it  Let $L$ be a lattice.
	Let $[x_0,x]\s L$  be a line-interval with atoms $x_i\
	(i\in I)$ (possibly $|I|=\infty$). Picking, for all $i\in I$,  any
	$p_i\in J(x_0,x_i)$ yields a line
	$l:=\{p_i|\ i\in I\}$
	with $\overline{l}=x$.  The map
	$p_i\mapsto x_0+p_i$ is a bijection from $l$ onto the atoms of
	$[x_0,x]$.}

\smallskip\rm
{\it Proof.}  Assuming $p_i\le x_j\ (i\not=j)$ yields the contradiction
$p_i\le x_ix_j=x_0$. Thus, since  the $x_i$'s are {\it all} lower covers
of $x$, one has $p_i+p_j=x\ (i\not=j)$. For all $p\in J(x)-l$ one has
$p\le x_i$ for some $i$, so $p+p_i\le x_i\not= x$. Hence $l$ is maximal with the property that any $p_i\neq p_j$ in $l$ yield the same join,
i.e. by definition $l$ is a line. $\square$

We say \underbar{$M_n$-line} for a line $\ell$ that derives from a line-interval $[x_0,x]$ in the way described in Lemma 7.1.
Because of $\ol{\ell}=x$ we also call $x$ a \underbar{line-top}.
 For instance $\ell_1=\{2,14,15,16\}$ in $L_1$ is of this sort, but not  line $\{2,6,8\}$ in $SM_{10}$. 
Any distinct join irreducibles $p_i,p_j$ of a $M_n$-line are perspective, since for any third  $p_k$ it holds (using the notation of Lemma 7.1) that
$(p_{i*},p_i)\nearrow (x_0,x_i)\nearrow (x_k,x)\nwarrow (x_0,x_j)\nwarrow (p_{j*},p_j).$

It would be convenient if in {\it  modular} lattices all lines were $M_n$-lines. 
The necessary condition of perspectivity is fulfilled:

\begin{enumerate}
	\item [(2)] {\it In each modular lattice any distinct join irreducibles $p,q$ of a line $\ell$ are perspective.}
\end{enumerate}

{\it Proof of (2).} Any distinct $p,q,r\in\ell$ satisfy $p+q=p+r=q+r=\ol{\ell}$.
 Since $q$ is 'weakly' prime (Ex. 2M) we have
$q\not\le q_*+r$. Therefore $q_*+r\prec q+r=\ol{\ell}$, and so $[q_*,q]\nearrow [q_*+r,\ol{\ell}]$. 
If we had $p\le q_*+r$, this would give the contradiction $q\le p+r\le q_*+r$.  Hence 
$p\not\le q_*+r$. As above this implies $[p_*,p]\nearrow [q_*+r,\ol{\ell}]$. This proves (2).

Is $\ell$ in (2) necessarily a $M_n$-line? This is not obvious since the join irreducibles of the line $\{2,6,8\}$ in $SM_{10}$ are also mutually perspective (a common upper transpose being $(9,1)$), notwithstanding the fact that $\{2,6,8\}$ is {\it no} $M_n$-line. Fortunately, as in (2), modularity again saves the day:

\medskip
{\bf Lemma 7.2 [HW,5.2]: }{\it  Let $L$ be a  modular lattice  and $p\ne q$  perspective join irreducibles.
	Then $[x_0,x]:=[p_*+q_*,p+q]$ is a line-interval and $x_0+p,\ x_0+q$ two of its
	atoms.}

{\it Proof (see Figure 7.1)}. Let $[e,f]$ be a common upper tranpose of $[p_*,p]$ and $[q_*,q]$.
For starters, $x:=p+q$ has at least two lower covers since $x\not\in J(L)$.
Let $x_0$ be the meet of {\it all} lower covers of $x$. Then $[x_0,x]$ is
a complemented (Thm. 2.2) modular lattice with
and $\delta([x_0,x])\ge 2$.
From $e\prec e+x=f$ follows
$x_0<ex\prec x$. From $p\not\le e$ follows $p\not\le x_0$, and so $\delta([px_0,p])\ge 1$. The isomorphy
$[x_0p,p]\simeq [x_0,x_0+p]$ implies that $x_0+p$ is join irreducible in
$[x_0,x_0+p]$, hence join irreducible in $[x_0,x]$,
hence an atom in the complemented interval $[x_0,x]$. Similarly $x_0+q$ is an atom
of $[x_0,x]$. From $(x_0+p)+(x_0+q)=x$ follows $\delta([x_0,x])=2$. Since $ex$ is a
third atom, $[x_0,x]$ is a line-interval. From $x_0\prec x_0+p$ follows $px_0=p_*$. Likewise $qx_0=q_*$, and so
$p_*+q_*\le x_0$.
From (6) in Section 2 follows $\delta(p+q)-\delta(p_*+q_*)\le 2$, and so $p_*+q_*=x_0$. $\square$

\begin{center}
	\includegraphics[scale=0.9]{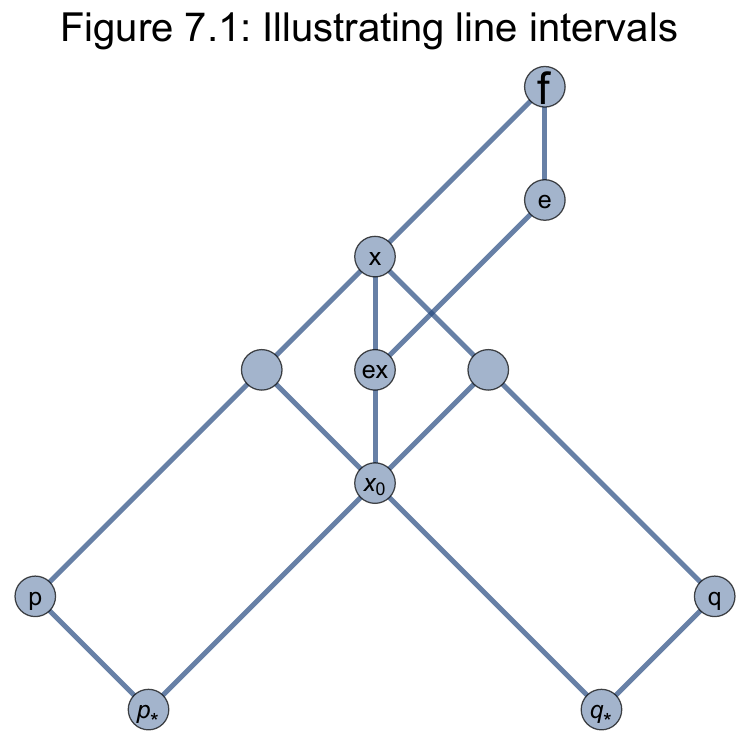} 
\end{center}

In view of (2) and Lemma 7.2 for each line $\ell$ in a  modular lattice we may put $\ul{\ell}:=x_0$, where $[x_0,x]$ is the line-interval coupled   to $\ell$.  
Since in the modular case, to which we stick from now on, each line is a $M_n$-line, we henceforth simply speak of 'lines'.

\begin{enumerate}
	\item [(3)] {\it In a modular lattice any distinct join irreducibles $p,q$ on a line $\ell$ satisfy $p_*+q_*=\ul{\ell}$.}
\end{enumerate}

This is immediate because $p,q$ in Lemma 7.2 can be {\it any} join irreducibles on the line $\ell$ that they determine, and we proved that  $p_*+q_*=x_0=\ul{\ell}$.

\smallskip
{\bf 7.2} Let $L$ be modular. Choosing
exactly one line for each line-interval yields a
family
$\Lambda\s \Pow(J(L))$.
We call the partial linear space
$\B:=(J(L),\Lambda)$  a \underbar{base of lines} of $L$.  For any $a\in L$ put $\B(a):=(J(a),\LL(a))$
where $\LL(a):=\{l\in\LL|\ \overline{l}\le a\}$. Obviously $\B(a)$ is a base of
lines for the lattice $[0,a]$. It is called
the base of lines  \underbar{induced} by $\B$.
\smallskip

\smallskip

{\bf 7.2.1} Return to $L_1$ of Figure 3.1. One has
$J(L_1)=\{2,3,4,6,7,10,12,14,15,16\}$ and $L_1$ has three line-intervals. The line
$l_1:=\{2,6,7\}$ corresponds to the line-interval $[3,9]$. There are two
lines $l_2':=\{6,14,15,16\}$ and $l_2:=\{2,14,15,16\}$ available for the line
interval $[10,18]$, we pick $l_2$. Finally, taking $l_3:=\{4,10,12\}$ for
$[9,17]$, one obtains a base of lines $\B_1:=(J(L),\LL)$ with $\LL:=\{l_1,
l_2,l_3\}$. 
What is the relevance of the connected components of the base of lines $\B_1$ of $L_1$? Take say
 $(p_*,p)=(3,7) $ and $(q_*,q)=(10,16)$. The fact that the lines containing $p$ and $q$ respectively, intersect in $2$, guides our choice of transpositions:

$(3,7)\nearrow (6,9)\searrow (0,2)\nearrow (14,18)\searrow (10,16).$

From this it should be clear that any two join irreducibles $p,q$ in a connected component of a base of lines of a modular lattice are projective. 

{\bf 7.3} Before we prove the converse in Theorem 7.3, let us investigate the situation in a {\it geometric} modular lattice $L$. Fix $x\in L\setminus\{0\}$. Since $[0,x]$ is atomistic, it follows from $(c)\Leftrightarrow (f)$ in Thm. 2.2 that $[0,x]$ is a line-interval iff $\delta(x)=2$ and $|[0,x]|\ge 5$. Each such line-interval $[0,x]$ houses exactly one line $\ell_x=At(x)$. Accordingly $L$ has a {\it unique} base of lines $\B_L=(At(L),\Lambda_L)$, where $\Lambda_L:=\{\ell_x:\ x\in L,\ \delta(x)=2\}$. If $p,q$ are in the same $\B_L$-component, then $p$ and $q$ are projective as illustrated in 7.2.1. Convesely, let $p,q\in J(L)$ be projective. Then $p,q$ are perspective by (8) in Section 4. By Lemma 4.2(b) they belong to a line-interval. So $p,q$ are $\B_L$-connected by virtue of a {\it single} line. To summarize, the $\B_L$-connected components are exactly the projectivity classes of join-irreducibles.
 
{\bf Theorem 7.3 $ [HW,2.6]$: }{\it Let $L$ be a modular lattice and $\B=(J(L),\LL)$
	a base of lines. 
	\begin{enumerate}
		\item [(a)] Then  $p,q\in J(L)$ are in the same connected component
		of $\B$ iff $(p_*,p)\approx(q_*,q)$.
		\item[(b)] The connected components
		$\B_i=(J_i,\LL_i)$ of $\B$ are isomorphic (as PLSes) to bases of lines for the $s(L)$
		congruence-simple factors $L/\theta_i$ of $L$.
	\end{enumerate}
	 }

\medskip

{\it Proof.} (a) One direction having been shown above, let $p',q'\in J(L)$ be arbitrary with $(p'_*,p')\approx(q'_*,q')$. We must show that $p',q'$ are $\B$-connected.
By join irreducibility $(p'_*,p'),(q'_*,q')$ cannot transpose down, and so the assumed projectivity has the form
$$(p'_*,p')\nearrow  (e_1,f_1)\searrow  (e_2,f_2)\nearrow \cdots\searrow  (e_{s-1},f_{s-1})\nearrow  (e_s,f_s)\searrow 
(q'_*,q').$$

We may assume (see Ex. 3A(i)) that all $f_{2i}$ are join irreducible. Therefore it suffices to
consider a two step projectivity $(p'_*,p')\nearrow  (e,f)\searrow  (q'_*,q')$.
Now $[x_0,x]:=[p'_*+q'_*,p'+q']$ is a line-interval (Lemma 7.2) which, among others, has the atoms $x_1:=x_0+p'$ and $x_2:=x_0+q'$.
By definition of a base of lines there is a $l\in\Lambda$ with
$[\underline{l},\overline{l}]=[x_0,x]$. Pick the points
$p,q\in l$ with $p\in J(x_0,x_1), q\in J(x_0,x_2)$. Because $(p_*,p)\approx (p_*', p')$ (they have $(x_0,x_1)$ as common upper transpose)  and $\delta([0,x_1])<\delta(L)$, by induction
$p$ and $p'$ are $\B(x_1)$-connected. Similarly $q$ and $q'$ are
$\B(x_2)$-connected. From $p,q\in l\in\LL$ follows that $p',q'$ are $\B$-connected.

As to (b), fix $i$ and let $f: L\to L/\theta_i$ be the canonical epimorphism. By (10) in Section 3, $f$ bijectivily
maps $J_i$ upon $J(L/\theta_i)$. Thus if $\ell\in\Lambda_i$ and $p,q\in\ell$ are distinct then $f(p),f(q)\in J(L/\theta_i)$ are distinct and $f(p)+f(q)=f(\ol{\ell})$. It follows that $f(\ol{\ell})$ is a line-top of $L/\theta_i$. Can $f(\ell)$ fail to be a line with line-top $f(\ol{\ell})$ because it isn't maximal w.r.t. the same-join-property? If so, extend $f(\ell)$ to a line $\ell'$ of $L/\theta_i$. But then, by Ex.3G, $\sigma(\ell')$ was a line of $L$ that properly contains $\ell$, a contradiction.
(Here $\sigma:L/\theta_i\to L$ is the smallest pre-images map.)
 Hence $f(\ell)$ is a line of $L/\theta_i$ for all $\ell\in\Lambda_i$. We see that $(J(L/\theta_i),\{f(\ell):\ \ell\in\Lambda_i\})$ can only fail to be a BOL of $L/\theta_i$ if there is some line-top $z\in L/\theta_i$ {\it different} from all line-tops $f(\ol{\ell})$. But then (why?) $\sigma(z)$ was a line-top of $L$ different from all line-tops $\ol{\ell}\ (\ell\in\Lambda_i)$, a contradiction.
 $\square$

{\bf 7.3.1}  Figure 7.2 renders a base of lines of the lattice $L_1$ in Figure 3.1(A). Since, say,  12 and 14 are in different connected components, it is by Thm. 7.3(a) impossible that the prime quotients $(12_*,12)$ and $(14_*,14)$ are projective. 
 By Thm. 7.3(b) the component $\{2,6,7,14,15,16\}$ in Figure 7.2 matches a BOL of some c-simple factor lattice of $L_1$.  Viewing that $\{2,6,7,14,15,16\}$ is $J(n)$ in Section 3.5, it follows that mentioned factor lattice is $L':=L_1/\theta_n$ depicted in Figure 3.1(B). Specifically, if $f$ is the natural epimorphism, then 
 
$\Big(f(J(n)),\ \{\ \{f(2),f(6),f(7)\},\ \{f(2),f(14),f(15),f(16)\} \} \Big)$

 must be a BOL of $L'$. Indeed, it is the BOL $(J(L'),\{\ \{p_1,p_2,p_3\},\{p_1,p_4,p_5,p_6\}\ \})$.

\begin{center}
	\includegraphics[scale=0.9]{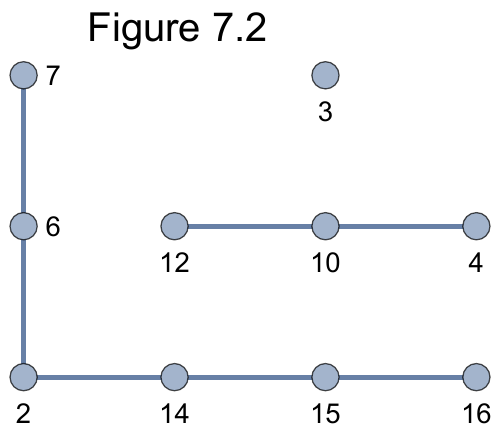} 
\end{center}

Exercise 7A asks to give a short proof of part (a) in Lemma 7.4 for the case of a {\it geometric} modular lattice.

{\bf Lemma 7.4: } {\it Let $L$ be modular with a base of lines $\B=(J(L),\LL)$.
\begin{enumerate}
	\item [(a)] Let $a\in L$  and $q,r\in J(L)$ incomparable such that
	$r\in J(a,  a+q)$. Then there is $p\in J(a)$ with $p+q=r+q$.
	\item[(b)] Let $S$ be any $\LL$-closed order ideal of $(J(L),\le)$. If $p,q\in S$ and $r\in J(L)$ and $r\le p+q$, then $r\in S$.
\end{enumerate}
}

{\it Proof.} As to (a), this is proven in [FH,Thm.4.2]. We like to give a more visual argument. In [G,p.85] the free modular lattice ('Dedekind's lattice') $FM(3)$ on $x,y,z$ is depicted. We replace $x,y,z$ by $a,r,q$ and recall that by assumption $r\le a+q$. This is equivalent to $a+q=a+r+q$. In the relabeled picture of $FM(3)$ the topmost prime quotient is $(a+q,a+r+q)$. Collapsing its 6-element projectivity class yields the lattice $S$ in Figure 7.3, which hence is (up to isomorphism) the largest sublattice of $L$ that $a,q,r\in L$ can possibly generate.
If the generated sublattice $\langle a,q,r\rangle$ is isomorphic to $S$, then picking any $p\in J(b,c)$ one sees that $p+q=r+q$. If $\langle a,q,r\rangle$ is smaller, then it must be (by universal algebra [KNT]) an epimorphic image of $S$, i.e $\langle a,q,r\rangle$ is obtained from $S$ by collapsing projectivity classes of prime quotients. One checks (Ex.7C) that also in these cases the sought $p$ exists.

\begin{center}
	\includegraphics[scale=0.7]{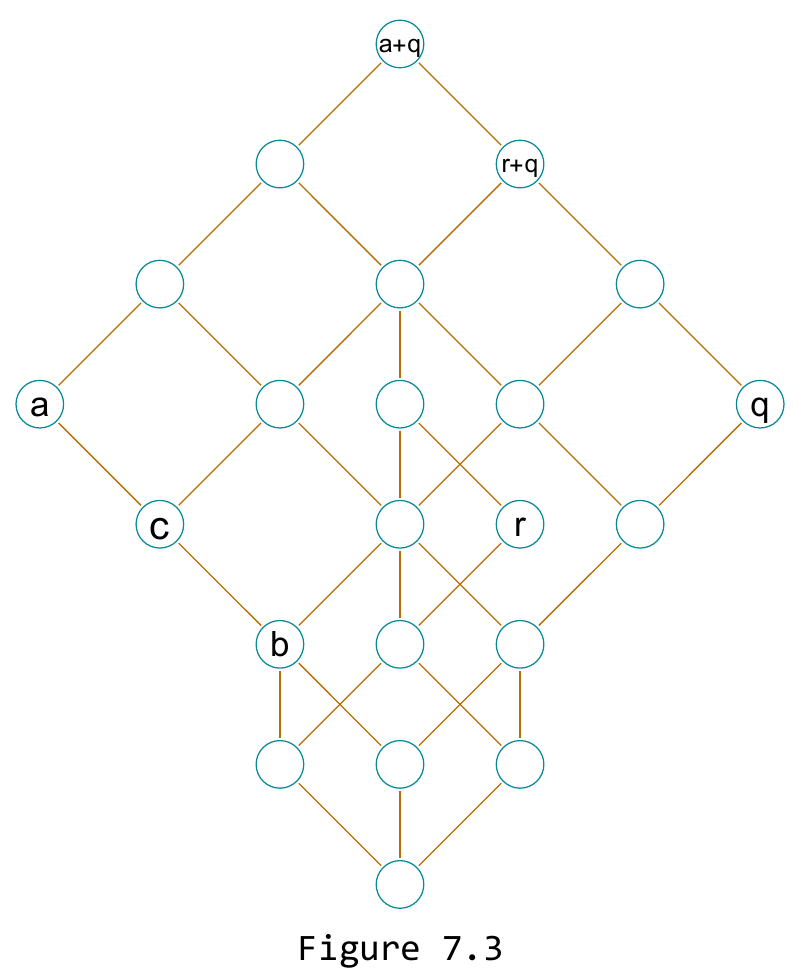} 
\end{center}

To prove (b) we use induction on the sum of the ranks of $p+q,p+r,q+r$.
{\it First case:} $p+r<p+q$ or $r+q<p+q$, say the latter. By (a) and
since $S$ is an order ideal, there is a $p'\in J(p)\s S$ with
$p'+q=r+q$. Since $\delta(p'+q)<\delta(p+q)$, induction applied to
$p',q,r$ yields $r\in S$. {\it Second case:} $p+r=p+q=r+q=:x$. Then 
$p,q,r$ can be extended to a line $l$ with $\overline{l}=x$. By definition of 'base of lines'
there is a $l'\in\Lambda$ that belongs to the line-interval
$[x_0,x]$. Hence there are points $p',q',r'\in l'$
belonging to the same atoms (say $x_1,x_2,x_3$) of $[x_0,x]$ as $p,q,r$. Induction\footnote{Notice that $p+q=p+q$ and $q+p'=q+r$, but $p+p'=x_1$ has smaller rank than $p+r=x$.} applied to
$p,q,p'$ (similarly $p,q,q'$) yields $p',q'\in S$. Thus $r'\in S$ since $S$ is
$\Lambda$-closed. Induction applied to $p',r',r$ finally yields $r\in S$. $\square$

\medskip

Theorem 7.5 generalizes a well known fact from projective geometry (recall Sec. 5.6.2). Benson and Conway [BC] proved the variant of Theorem 7.5 where $\LL=\LL_L$ is the set of {\it all} lines of $L$. This allows to give a  shorter\footnote{This will make a guided exercise in some later version of Part B.} proof. For any BOL $\B=(J(L),\Lambda)$ of $L$ we denote by $\La(\B,\le)$ the closure system of all $\Lambda$-closed order ideals of $(J(L),\le)$.

{\bf Theorem 7.5 [HW,2.5]: }{\it Let $\B=(J(L),\Lambda)$ be a base of lines of the 
 modular lattice $L$. Then $a\mapsto J(a)$ is a lattice isomorphism
from $L$ onto $\La(J(L),\LL,\le)$.}

\medskip

{\it Proof of Theorem 7.5.}
Finite length guarantees $a=\sum J(a)$ for all $a\in L$.
Each set $J(a)$ is an order ideal of $(J(L),\le)$  which is $\LL$-closed because from $\ell\cap J(a)\supseteq\{p,q\}$ and $r\in\ell$ follows $r\le p+q\le a$, i.e. $r\in J(a)$.
Therefore $f(a):=J(a)$ is a well defined map $L\to \La(\B,\le)$. It is an order embedding in view of $a\le b\Leftrightarrow J(a)\s J(b)$. If we manage to show $f$ is
 {\it onto}, $f$ will be a lattice isomorphism (see Ex.2C).

 So let $S\in \La(\B,\le)$
be arbitrary. Trivially $S\s J(\sum S)$. We tackle the proof of the
other inclusion by choosing
a finite subset  $\{p_1,\dots,p_m\}\s S$
with $p_1+\cdots+p_m=\sum S$ and by inducting on $k$ to show that
$J(p_1+\cdots p_k)\s S$ for all $1\le k\le m$.
As to $k=1$, from $p_1\in S$ follows $J(p_1)\s S$ since $S$ is an order ideal of $(J(L),\le)$.
 For the induction step, put $a:=p_1+\cdots p_k,\ q:=p_{k+1}$, and assume that $J(a)\s S$. Fix any $r\in J(a+q)$. We need to show that $r\in S$.
 {\it Case 1:} $r> q$. Then from $a+r_*\ge a+q\ge r$ and the weak primality of $r$ (Ex.2M) follows the contradiction $a\ge r$.
 {\it Case 2:} $r\le q$. Then $r\in J(q)\s S$. {\it Case 3:} $r,q$ are incomparable. Then by Lemma 7.4(a) there is $p\in J(a)$ with $p+q=r+q$. From $r\le p+q$ and $p,q\in S$ follows $r\in S$ by Lemma 7.4(b). $\square$
 
{\bf 7.4}  Here comes a summary of how geometric, respectively modular lattices decompose according to the behaviour of their join irreducibles:
 
 \begin{itemize}
 	\item For {\bf both} types of lattices $L$ the projectivity classes of join irreducibles match the congruence-simple (c-simple) factors of $L$.
 	\item For {\bf geometric} lattices  "c-simple $\LRa$ directly irreducible", and for any two atoms it holds that "projective $\LRa$ perspective". 
 	 The c-simple factors match the connected components of the  matroid that lives on the atoms.
 	\item For {\bf modular} lattices the c-simple factors match the $\B$-connected components w.r.t. any base of lines BOL $\B$. Furthermore, "c-simple $\LRa$ subdirectly irreducible".
 	\item For {\bf geometric and modular} lattices the c-simple factors match projective spaces. For any two atoms it holds that "projective $\LRa$ perspective $\LRa$ collinear". 
 \end{itemize}
 \bigskip
 
 {\bf 7.5. }  Recall that by Theorem 7.5 each modular lattice $L$ is isomorphic to the closure system $\La(\B,\le)$ of all $\LL$-closed order ideals of $(J(L),\le)$. Here $\B=(J(L),\LL)$ is any BOL of $L$. A natural implicational base (Sec. 2.3) for $\La(\B,\le)$  is
 
 $(4)\quad \Sigma_{nat}:=\{\{p\}\to J(p):\ p\in J^*\}\cup \{\{p,q\}\to\ell:\ \ell\in\LL,\ p\neq q\ in\ \ell\}.$ 
 \bigskip
 
 As to the definition of $J^*$, when $p$ is an atom, the implication $\{p\}\to J(p)$ boils down to $\{p\}\to\{p\}$ and whence can be dropped.
 Hence it suffices that $p$ ranges over $J^*:= J(L)\setminus At(L)$. It turns out that $\Sigma_{nat}$ is a \underbar{minimum} implicational base in the sense that $|\Sigma|\ge|\Sigma_{nat}|$ for each implicational base 
 $\Sigma$ of $\La(\B)$. 
 
 Generally, the \underbar{size} of an implicational base $\Sigma=\{A_1\to B_1,\ldots,A_t\to B_t\}$ of any closure system is defined as
 $s(\Sigma):= |A_1|+\cdots+|A_t|+|B_1|+\cdots+|B_t|.$
 An implicational base $\Sigma'$ of some fixed closure system $\La$ is \underbar{optimal} if $s(\Sigma)\ge s(\Sigma')$ for all implicational bases $\Sigma$ of $\La$.
 
 Returning to $\La(\B,\le)$, while $\Sigma_{nat}$ is minimum, it is far from being optimal. For instance, if $L$ is the lattice in Figure 3.1, then 
 
 $\{4\}\to\{2\},\ \{6\}\to\{3\},\ \{7\}\to\{3\},\ \{10\}\to\{7\},\ \{12\}\to\{2,6,7\},\ \{14\}\to\{10\},\\ \{15\}\to\{10\},\ \{16\}\to\{10\},\
 \{4,10\}\to\{12\},\ \{4,12\}\to\{10\},\ \{10,12\}\to\{4\},\\ \{2,14\}\to\{15\},\ \{2,15\}\to\{14\},\ \{2,16\}\to\{14\},\\ \{14,15\}\to\{16\},\ \{14,16\}\to\{15\},\ \{15,16\}\to\{2\},$
 
 is one of many optimal implicational bases, all of which  having $s(\Sigma_{opt})=47$. For details, the reader is referred to [W5].
 
 \bigskip

{\bf 7.6 Exercises}

{\bf Exercise 7A.} Give a short proof of Lemma 7.4(a) for geometric modular lattices $L$.

{\bf Solution of Ex. 7A.} Since $q\in At(L)$ we have $a\prec a+q$. Hence $r\in J(a,a+q)$ forces $a+q=a+r$. Therefore $(r_*,r)\nearrow (a,a+q)\nwarrow (q_*,q)$, i.e. $r,q$ are perspective. By Lemma 4.2(b), $x:=r+q$ yields a line-interval $[0,x]$. From $x\prec a+x\ (=a+q)$ follows $ax\prec x$. Hence $p:=ax$ is an atom below $a$ with $p+q=r+q$.

{\bf Exercise 7B.} Find bases of lines for all lattices in Figure 3.2 of Exercise 3E.

{\bf Exercise 7C.} Complete the proof of Lemma 7.4(a). This concerns epimorphic images of Figure 7.3.

\bigskip

\centerline{\bf 8. Calculating the submodule lattice of a finite $R$-module}

Let us consider any poset and any PLS that share\footnote{Of course BOL's of modular lattices fit that hat, but for the time being we adopt a more abstract setting.} a common universe, thus $(E,\le)$ and $(E,\LL)$.
For $a\in E$ we put $a\downarrow:=\{b\in E:\ b\le a\}$.
 We aim for a compressed representation of the closure system $\La(E,\le,\LL)$ of all $\LL$-closed order ideals. Since the members of
$\La(E,\le,\LL)$ are exactly the sets $X\s E$ which are $\Sigma$-closed w.r.t. 

$\Sigma=\{a\ra (a\downarrow): a\in E\}\cup\{\{a,b\}\ra\ell:\ a,b\in\ell\in\LL,a\neq b\},$

we can obtain the desired representation of $\La(E,\le,\LL)$ in polynomial total time by applying the Horn $n$-algorithm of [W7] to $\Sigma$. However, here the special types of implications invite shortcuts. First, the implications induced by the partial ordering of $E$ have {\it singleton premises}. Also for each fixed $\ell$  the ${|\ell|}\choose{2}$ implications $\{a,b\}\to\ell$ can likely be trimmed since they all have the {\it same conclusion} $\ell$.

{\bf 8.1} For instance, let
$E:=\{p_1,\ldots,p_7\}$ and let $(E,\le)$ and $(E,\LL)$ be as in Figure 8.1.

\begin{center}
	\includegraphics[scale=0.85]{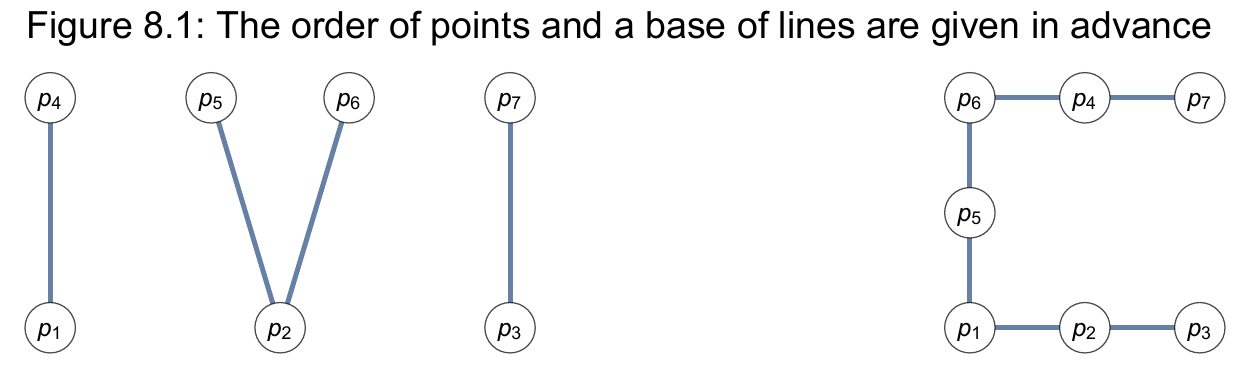} 
\end{center}

The preliminary system $\La(E,\le)$ of all order ideals (viewed as bitstrings $(x_1,...,x_7)$) can be rendered as a disjoint union of so called $(a,b)$-{\it rows} $r_i$. Details being provided in [W6] we will be brief. We split $\La(E,\le)$ as $r_1\uplus r_2$ where $r_1$ contains the order ideals containing $p_2$ (i.e. $x_2=1$), and $r_2$ the ones with $x_2=0$. By definition of 'order ideal' $x_2=0$ forces $x_5=x_6=0$, whereas $x_2=1$ allows for $x_5=x_6=2$, where '2' is a don't-care symbol (freely replaceable by 0 or 1).
Furthermore, these conditions must be satisfied by each wannabe order ideal:

$ (x_4=1\Ra x_1=1)\ and\ (x_7=1\Ra x_3=1)$

There are three 0,1-assignments to $x_4,x_1$ that satisfy $(x_4=1\Ra x_1=1)$, namely
$(x_4,x_1)=(1,1),(0,0),(0,1)$.We abbreviate this by the wildcard $(a_1,b_1)$ that appears in both $r_1$ and $r_2$. Likewise there are three 0,1-assignments for $(x_7,x_3)$, which explains the wildcard $(a_2,b_2)$ in $r_1$ and $r_2$. Since the bits $x_5,x_6$ in $r_1$ are not restricted, we see that $|r_1|=3\cdot 3\cdot 4=36$. 
Similarly $|r_2|=9$, and so $|\La(E,\le)|=36+9=45$.

\begin{tabular}{c|c|c|c|c|c|c|c|c  }
	& $p_1$ & $p_2$ & $p_3$& $p_4$ & $p_5$ & $p_6$ &$p_7$ &  \\ \hline 
	&  &  &  & &  &   &  & \\ \hline
	$r_1=$ & $b_1$ & $\bf 1$ &$b_2$ &$a_1$ & $2$ & $2$ & $a_2$ & 
	{\it pending} $\ell_1$ \\ \hline
	$r_2=$ & $b_1$ & $\bf 0$ & $b_2$ & $a_1$ &$0$& $0$ &$a_2$ &
	{\it pending} $\ell_1$   \\ \hline
	
	&  &  &  & &  &   &  & \\ \hline
	$r_3=$ & ${\bf 0}$ & $1$ & ${\bf 0}$ & $0$ & $2$ & $2$ & $0$ & 
	{\it pending} $\ell_2$ \\ \hline
	$r_4=$ & ${\bf 1}$ & $1$ & ${\bf 1}$ & $2$ &$2$& $2$ & $2$ &
	{\it pending} $\ell_2$   \\ \hline
	$r_2=$ & $b_1$ & $0$ & $b_2$ & $a_1$ &$0$& $0$ &$a_2$ &
	{\it pending} $\ell_1$   \\ \hline
	
	&  &  &  & &  &   &  & \\ \hline
	$r_5=$ & $0$ & $1$ & $0$ & $0$ & $\epsilon$ & $\epsilon$ & $0$ & 
	{\it final}  \\ \hline
	$r_4=$ & $1$ & $1$ & $1$ & $2$ &$2$& $2$ & $2$ &
	{\it pending} $\ell_2$   \\ \hline
	$r_2=$ & $b_1$ & $0$ & $b_2$ & $a_1$ &$0$& $0$ &$a_2$ &
	{\it pending} $\ell_1$   \\ \hline
	
	&  &  &  & &  &   &  & \\ \hline
	$r_6=$ & $1$ & $1$ & $1$ & $2$ & ${\bf 0}$ & ${\bf 0}$ & $2$ & 
	{\it pending} $\ell_3$  \\ \hline
	$r_7=$ & $1$ & $1$ & $1$ & $2$ &${\bf 1}$& ${\bf 1}$ & $2$ &
	{\it pending} $\ell_3$   \\ \hline
	$r_2=$ & $b_1$ & $0$ & $b_2$ & $a_1$ &$0$& $0$ &$a_2$ &
	{\it pending} $\ell_1$   \\ \hline
	
	&  &  &  & &  &   &  & \\ \hline
	$r_8=$ & $1$ & $1$ & $1$ & $\epsilon$ & $0$ & $0$ & $\epsilon$ & 
	{\it final}   \\ \hline
	$r_7=$ & $1$ & $1$ & $1$ & $2$ &${\bf 1}$& ${\bf 1}$ & $2$ &
	{\it pending} $\ell_3$   \\ \hline
	$r_2=$ & $b_1$ & $0$ & $b_2$ & $a_1$ &$0$& $0$ &$a_2$ &
	{\it pending} $\ell_1$   \\ \hline
	
	&  &  &  & &  &   &  & \\ \hline
	$r_9=$ & $1$ & $1$ & $1$ & $d$ &$1$& $1$ & $d$ &
	{\it final}    \\ \hline
	$r_2=$ & $b_1$ & $0$ & $b_2$ & $a_1$ &$0$& $0$ &$a_2$ &
	{\it pending} $\ell_1$   \\ \hline
	
	&  &  &  & &  &   &  & \\ \hline
	$r_{10}=$ & ${\bf 0}$ & $0$ & ${\bf 0}$ & $0$ & $0$ & $0$ & $0$ & 
	{\it final}   \\ \hline
	$r_{11}=$ & ${\bf 0}$ & $0$ & ${\bf 1}$ & $0$ &$0$& $0$ & $2$ &
	{\it final}    \\ \hline
	$r_{12}=$ & ${\bf 1}$ & $0$ & ${\bf 0}$ & $2$ &$0$& $0$ &$0$ &
	{\it final}    \\ \hline
	
\end{tabular}
\bigskip

{\sl Table 8.1: Steps towards enhancing the (a,b)-algorithm}

In order to sieve the $\LL$-closed order ideals from $r_1\uplus r_2$, we "impose" the lines $\ell_1=\{p_1,p_2,p_3\},\ell_2=\{p_1,p_5,p_6\}$ and $\ell_3=\{p_4,p_6,p_7\}$ one after the other. Those bitstrings $x\in r_1$ that are $\{\ell_1\}$-closed (in the sense of Sec. 5.5) must either satisfy 

(a) $x_1=x_3=0$, or (b) $x_1=x_3=1$.

 But this is not enough. In case (a) it follows from $b_1=0$ that $a_1=0$; and likewise $b_2=0$ forces $a_2=0$. Hence the type (a) bitstrings in $r_1$ that are $\{\ell_1\}$-closed are collected in $r_3$. Likewise the type (b) bitstrings in $r_1$ that are $\{\ell_1\}$-closed are collected in $r_4$. Our {\it working stack} now contains, top down, the three rows $r_3,r_4,r_2$. In both $r_3,r_4$ the line $\ell_2$ is {\it pending} to be imposed, in $r_1$ it is still $\ell_1$.
According to the last-in-first-out (LIFO) principle we always process the top row of the working stack, i.e. now $r_3$. The bitstrings $x\in r_3$ which are $\{\ell_2\}$-closed are exactly the ones which do {\it not} have $x_5=x_6=1$. In other words, we need 'at most one 1 among $x_5,x_6$'. We indicate that with the wildcard\footnote{The old-school alternative would replace $r_5$ by three ordinary bitstrings.} $(\epsilon,\epsilon)$. The so obtained row $r_5$ is $\{\ell_2\}$-closed by construction. It happens to be $\{\ell_3\}$-closed as well, and so is {\it final} in the sense that $r_5\s\La(\B,\le)$. We remove $r_5$ from the working stack and store it in a safe place. This frees $r_4$ to be processed. It is clear that the $\{\ell_2\}$-closed bitstrings
 $x\in r_4$ are exactly the ones in $r_6\uplus r_7$. Imposing $\ell_3$ on $r_6$ yields the final row $r_8$. Imposing $\ell_3$ on $r_7$ yields the final row $r_9$. We used yet another wildcard here. By definition $dd\cdots d$ means 'all 1's or all 0's'.
In order to impose $\ell_1$ on the last row $r_2$ in the working stack, we cannot avoid to spell out the three possibilities $(0,0),(0,1),(1,0)$ for $(b_1,b_2)$. This yields the final rows $r_{10},r_{11},r_{12}$, and our algorithm  terminates. 
We conclude that \smallskip

\centerline{$\La(E,\le,\LL)=r_5\uplus r_8\uplus r_9\uplus r_{10}\uplus r_{11}\uplus r_{12}$ and $|\La(E,\le,\LL)|=3+3+2+1+2+2=13.$}

(Generally such a a compressed representation  yields the cardinality at once.) The diagram of the lattice $L(E,\le,\LL)$ can be erected almost in 3D-printing fashion. Namely, starting with the zero bitstring (which is contained in $r_{10}$ and represents the empty order ideal), look for the minimal bitstrings in the remaining set system. They are $(1,0,0,0,0,0,0),(0,1,0,0,0,0,0),(0,0,1,0,0,0,0)$ and are contained in $r_{12},r_5,r_{11}$ respectively. The next layer to be 'printed' consists of bitstrings contained in $r_{12},r_5,r_8,r_5,r_{11}$ respectively. And so it goes on. It is easy to connect each bitstring with the correct bitstrings in the layer below. In this way one obtains the  lattice in Fig.8.2(A). We labeled each node (=$\LL$-closed order ideal) with the row $r_i$ that contained the bitstring  matching it. The labelling in Fig.8.2(B) shows that this lattice indeed has a poset of join irreducibles, and a base of lines $\LL$, as in Figure 8.1. Notice that $\ol{\ell_1}=z,\ \ol{\ell_2}=y,\ \ol{\ell_3}=x$.

\begin{center}
	\includegraphics[scale=0.7]{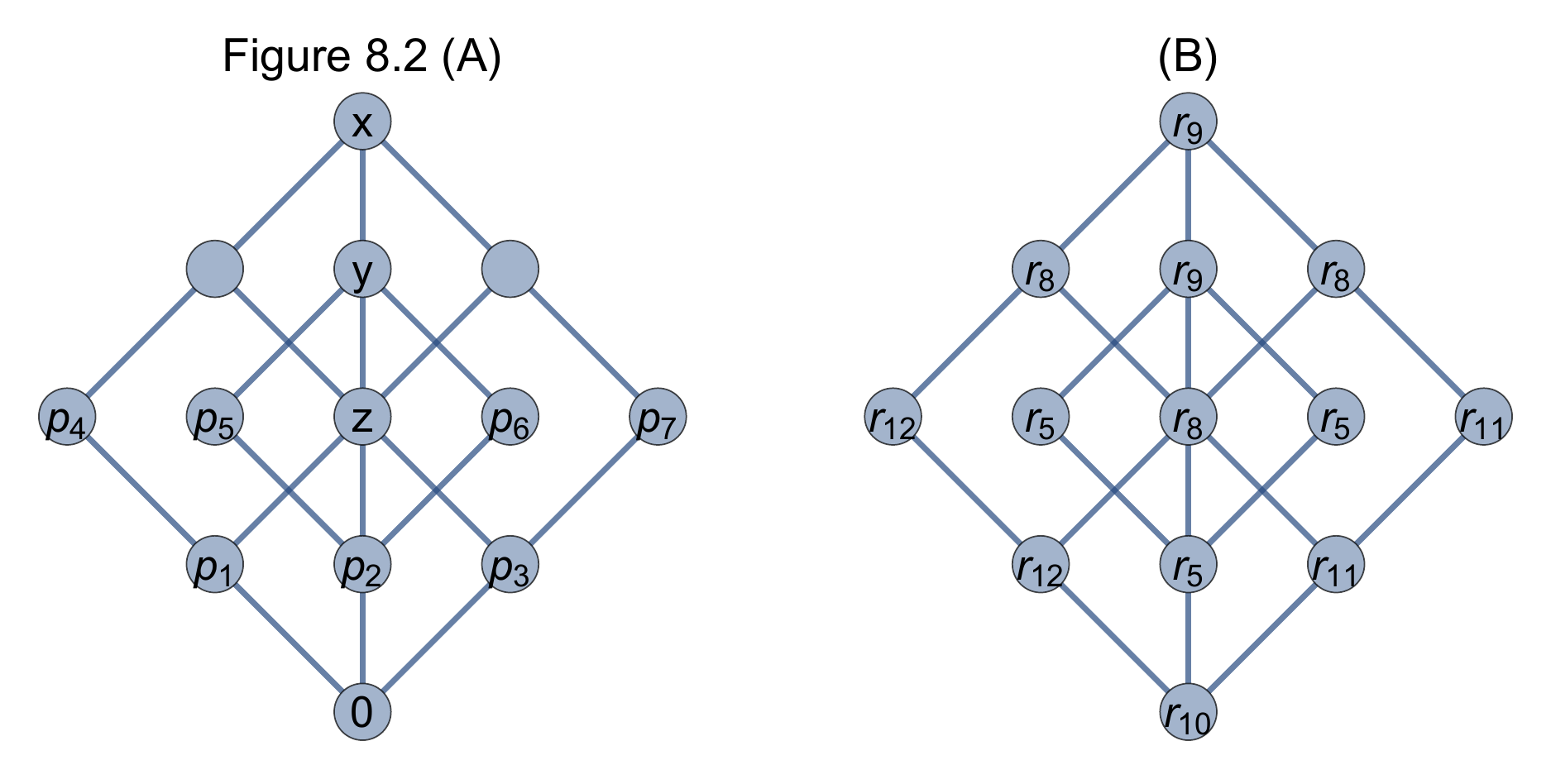} 
\end{center}

{\bf 8.2} Our toy example indicates that upon running the (a,b)-algorithm, in order to impose the lines of $\LL$, it is handy to employ $d$-wildcards and $\epsilon$-wildcards. Different from Table 8.1 it may happen that a line to be imposed on $r$ only touches upon the 2's of $r$. One can then employ the $\ell$-wildcard $(\ell,...,\ell)$ which by definition means "at most one 1, or all 1's". Upon embracing the $\ell$-wildcard one finds that the $g$-wildcard would be handy as well. By definition $(g,...,g)$ means "exactly one 1 here". Fortunately now the inflation of wildcards stops. Specifically, consider any $(0,1,2,a,b,d,\epsilon,g,\ell)$-row $r$ and any $\ell$-wildcard of length $\lambda$. Then, we claim,  the $\ell$-wildcard can be imposed on $r$ by substituting $r$ with at most $\lambda+2$ suitable, pairwise disjoint  $(0,1,2,a,b,d,\epsilon,g,\ell)$-rows.
Rather than giving a tedious formal proof, we illustrate the procedure, call it  $(0,1,2,a,b,d,\epsilon,g,\ell)$-{\it algorithm}, in Table 8.2.

\begin{tabular}{|c|c|c|c|c|c|c|c|c|c|c|c|c|c|c|c|c|c|c|c|c|c|c|  }
	1 & 2 & 3& 4 & 5  & 6 &7 &8 &9  &10 &11  &  &12 &13	&14  &15 &16  &17  &18 &19  &20 &21 &22 \\ \hline\hline
	
	$d$ & $d$ & $\epsilon$& $\ell_1$ & $\ell_1$ & $\ell_2$ &$\ell_2$ & $a$ & $b'$ & $g$& $2$ &  &
	$d$ &$\epsilon$	& $\epsilon$ &$\ell_1$ & $\ell_2$ & $\ell_2$ & $b$& $b$ &$a'$  & $g$ &$g$ \\ \hline
	
	&  &  & &  &   & &  &  & &  &   & &  &   & &  &  & &  &   & & \\ \hline
	
	$0$ & $0$ & $0$& $0$ & $0$ & $0$ &$0$ & $0$ & $0$ & $0$& $0$ &  &
	$0$ &$\epsilon$	& $\epsilon$ &$2$ & $\epsilon'$ & $\epsilon'$ & $2$& $2$ &$0$  & $g$ &$g$ \\ \hline
	
	&  &  & &  &   & &  &  & &  &   & &  &   & &  &  & &  &   & & \\ \hline
	
	$0$ & $0$ & ${\bf 1}$& $0$ & $0$ & $0$ &$0$ & $0$ & $0$ & $0$& $0$ &  &
	$0$ &${\bf 0}$	& ${\bf 0}$ &$2$ & $\epsilon'$ & $\epsilon'$ & $2$& $2$ &$0$  & $g$ &$g$ \\ \hline
	
	$0$ & $0$ & $0$& ${\bf g}$ & ${\bf g}$ & $0$ &$0$ & $0$ & $0$ & $0$& $0$ &  &
	$0$ &$\epsilon$	& $\epsilon$ &${\bf 0}$ & $\epsilon'$ & $\epsilon'$ & $2$& $2$ &$0$  & $g$ &$g$ \\ \hline
	
	$0$ & $0$ & $0$& $0$ & $0$ & ${\bf g}$ &${\bf g}$ & $0$ & $0$ & $0$& $0$ &  &
	$0$ &$\epsilon$	& $\epsilon$ &$2$ & ${\bf d}$ & ${\bf d}$ & $2$& $2$ &$0$  & $g$ &$g$ \\ \hline
	
	$0$ & $0$ & $0$& $0$ & $0$ & $0$ &$0$ & ${\bf 1}$ & $0$ & $0$& $0$ &  &
	$0$ &$\epsilon$	& $\epsilon$ &$2$ & $\epsilon'$ & $\epsilon'$ & ${\bf 1}$& ${\bf 1}$ &$0$  & $g$ &$g$ \\ \hline
	
	$0$ & $0$ & $0$& $0$ & $0$ & $0$ &$0$ & $0$ & ${\bf 1}$ & $0$& $0$ &  &
	$0$ &$\epsilon$	& $\epsilon$ &$2$ & $\epsilon'$ & $\epsilon'$ & $2$& $2$ &${\bf 2}$  & $g$ &$g$ \\ \hline
	
	$0$ & $0$ & $0$& $0$ & $0$ & $0$ &$0$ & $0$ & $0$ & ${\bf 1}$& $0$ &  &
	$0$ &$\epsilon$	& $\epsilon$ &$2$ & $\epsilon'$ & $\epsilon'$ & $2$& $2$ &$0$  & ${\bf 0}$ &${\bf 0}$ \\ \hline
	
	$0$ & $0$ & $0$& $0$ & $0$ & $0$ &$0$ & $0$ & $0$ & $0$& ${\bf 1}$ &  &
	$0$ &$\epsilon$	& $\epsilon$ &$2$ & $\epsilon'$ & $\epsilon'$ & $2$& $2$ &$0$  & $g$ &$g$ \\ \hline
	
	&  &  & &  &   & &  &  & &  &   & &  &   & &  &  & &  &   & & \\ \hline
	
	$1$ & $1$ & $1$& $1$ & $1$ & $1$ &$1$ & $1$ & $1$ & $1$& $1$ &  &
	$1$ &$0$	& $0$ &$1$ & $1$ & $1$ & $1$& $1$ &$2$  & $0$ &$0$ \\ \hline
		
\end{tabular}
\bigskip

{\sl Table 8.2 Splitting a typical $(0,1,2,a,b,\epsilon,d,g,\ell)$-row by imposing an $\ell$-wildcard\\ occupying positions 1 to 11.}

Here the top row $r:=(d,d,\epsilon,...,g)$ undergoes the imposition of a $\ell$-wildcard with position-set $\{1,2,..,11\}$. In a nutshell this is achieved as follows. All $X\in R$ with $X\cap\{1,...,11\}=\es$ can be collected in a single row of type $(0,...,0,*,..*)$, all $X\in r$ with $\{1,2,...,11\}\s X$ in a single row of type $(1,...,1,*,...,*)$,
but it takes seven (disjoint) rows to gather all $X\in r$ with $|X\cap\{1,2,...,11\}|=1$. We leave it to the reader to unravel the details.

{\bf 8.3} Let us  calculate the subgroup lattice $L_0$ of the Abelian group (=$\Z$-module)
$\Z_4\times\Z_4$. If we can find the poset $(J(L_0),\le)$ and any base of lines $\LL$, then the  $(0,1,2,a,b,d,\epsilon,g,\ell)$-algorithm can be applied. In each $R$-submodule lattice the join-irreducibles are the 1-generated submodules $\langle x\rangle$. In our case they are

$p_1=\langle(0,2)\rangle,\ p_2=\langle(2,2)\rangle,\
p_3=\langle(2,0)\rangle,\ p_4=\langle(0,1)\rangle,\ p_5=\langle(1,1)\rangle$

$p_6=\langle(1,3)\rangle\  (=\{(1,3),(2,2),(3,1),(0,0)\}),\  p_7=\langle(1,0)\rangle. $

These join-irreducibles happen to be  ordered as in Figure 8.2(A). We now indicate a naive method (which likely can be trimmed) to come up with a base of lines in a general  submodule lattice $L'$. First  we calculate all $M_n$-elements. Among all submodules of type $p_i+p_j\ (i\neq j)$ they are exactly the ones which admit a third $p_k$ such that\footnote{For general rings $R$ and $R$-modules $W$ we must postulate that the finitely many 1-generated submodules $\langle a\rangle$ are known sufficently well to efficiently 
	verify equalities of type 
	$\langle a\rangle+\langle b\rangle=\langle c\rangle+\langle d\rangle$.   
Of course that holds iff there are $\alpha_1,...,\delta_2\in R$ with $a=\gamma_1 c+\delta_1 d,\ b=\gamma_2 c+\delta_2 d,\ c=\alpha_1 a+\beta_1 b,\ d=\alpha_2 a+\beta_2 b$. The question is, how fast $\alpha_1,...,\delta_2\in R$ can be found.} $p_i+p_j=p_i+p_k=p_j+p_k$.
For each  $M_n$-element $x$ extend the already found partial line $\{p_i,p_j,p_k\}$ to a set $\ell_x\s J(L)$ maximal with the property that any two distinct elements of $\ell_x$ have join $x$. Then the family $\LL'$ of all these lines $\ell_x$ is a base of lines. In our particular $L_0$ there are three $M_3$-elements $x,y,z$ (see Fig.8.2(A)) and the family $\LL$ as in Fig.8.1 is one possible base of lines.

\bigskip

{\bf 8.4} It is tempting to generalize matters. But it depends on this

{\bf Open Question:} {\it  Given any submodules $X_1,..,X_t$ of a finite-length $R$-module $W$, is there an efficient method to calculate the join-irreducibles of the generated submodule lattice $L=\langle X_1,...,X_t\rangle$?}

If yes, then the $(0,1,2,a,b,\epsilon,d,g,\ell)$-algorithm can be applied to the join-irreducibles. It would be interesting to compare our method with the wholly different approach in [LMR].

{\bf 8.5}  It is worthwile pointing out that the Open Question becomes easy in the distributive case. Specifically, given sets $X_1,..,X_t\in{\cal P}(W)$, how can the generated (necessarily distributive) sublattice $D=\langle X_1,...,X_t\rangle$ be calculated in a compressed format?

By distributivity each $Y\in D$ is a union of intersections of the $X_i$'s, and so the join irreducibles $A\in D$ are {\it among} these intersections. 
Specifically, putting

$(1)\quad A_v:=\bigcap\{X_i:\ 1\le i\le t,\ v\in X_i\},$

we claim that

$(2)\quad J(D)=\{A_v:\ v\in W\}.$

Indeed, any intersection $\bigcap\{X_i:\ i\in I\}$ strictly contained in $A_v$ does not contain $v$ (because $v\not\in X_i$ for some $i$ by definition pf $A_v$). Hence $A_v$ isn't a union of such intersections, hence $A_v$ is join-irreducible. Conversely, let $A\in J(D)$ be arbitrary. We know that $A=\bigcap\{X_k:\ k\in K\}$ for some index set $K\s [t]$.  Furthermore, let $A_*$ be the unique lower cover of $A$ in $D$ and fix any $v\in A\setminus A_*$. From $v\in A\s X_k$ for all $k\in K$ follows that $A_v\s A$. Hence either $A_v=A$ or $A_v\s A_*$. The latter is impossible because $v\not\in A_*$. This proves (2).

Because  $J(D)$, partially ordered by inclusion, determines $D$ via $D\simeq \La(J(D),\s)$ 
(Birkhoff's Theorem), it suffices to give a compressed representation of $\La(J(D),\s)$ . This is achieved by the (a,b)-algorithm. See Exercises 8A and 8B.
 
\bigskip

{\bf 8.6 Exercises}

{\bf Exercise 8A:}  Put $W:=\{a,b,c,d,e,f,g,h,k\}$ and consider the subsets $X_i\s W$ defined by the rows of the Table below. Calculate the join-irreducibles of $D=\langle X_1,...,X_8\rangle$ as sketched in 8.5.

\begin{tabular}{c|c|c|c|c|c|c|c|c|c|  }
&a	& b & c & d& e & f & g &h & k \\ \hline 
	& & &  &  & &  &   &  & \\ \hline
$X_1=$  &1	 & 1 & 1 &1 &1 & 1 & 0 & 0 &1  \\ \hline
$X_2=$  &1	 & 1 & 1 &1 &1 & 1 & 0 & 1 &0  \\ \hline
$X_3=$  &1	 & 1 & 0 &1 &1 & 1 & 0 & 0 &0  \\ \hline
$X_4=$  &0	 & 1 & 1 &1 &0 & 1 & 1 & 0 &0  \\ \hline
$X_5=$  &1	 & 1 & 0 &0 &1 & 0 & 0 & 0 &0  \\ \hline
$X_6=$  &0	 & 1 & 0 &1 &0 & 1 & 0 & 0 &0  \\ \hline
$X_7=$  &1	 & 1 & 0 &1 &1 & 1 & 1 & 1 &0  \\ \hline
$X_8=$  &0	 & 0 & 0 &1 &0 & 1 & 1 & 0 &1  \\ \hline
	\end{tabular}

{\bf Exercise 8B:} Having computed $J(D)$ in Ex.8A, find a compressed enumeration of $D$ with the (a,b)-algorithm. (This may require to consult [W6].)

\vspace{2cm}

 \centerline{\bf 9. Bases of lines: Localization, and two types of cycles}

{\bf 9.1. } Let $\B$ be a base of lines of
 the modular lattice $L$ and $a\prec b$ a covering.
 The \underbar{localization} of $\B=(J(L),\LL)$ to $a\prec b$ is the
 space $\B(a,b):=(J(a,b),\LL(a,b))$, where $\LL(a,b):=\{l\cap J(a,b)|\ l\in\LL,
 \overline{l}\le b,\overline{l}\not\le a\}$. Consider a line $l':=l\cap
 J(a,b)$ of $\B(a,b)$. Because $a+\overline{l}=b$ implies $\underline{l}
 \prec \overline{l}a\prec\overline{l}$, there is a $p\in l$ with $p\le a$.
 Any other $q\in l$ with $q\le a$ yields the contradiction $\overline{l}=
 p+q\le a$. So $l=l'\cup\{p\}$. Other than the induced base of lines $\B(a)$
 the localization $\B(a,b)$ is not the base of lines
 of any lattice (in particular $\vv l'\vv=2$ is possible), but it will be an
 important gadget.
 \smallskip
 
  \begin{center}
 	\includegraphics[scale=0.8]{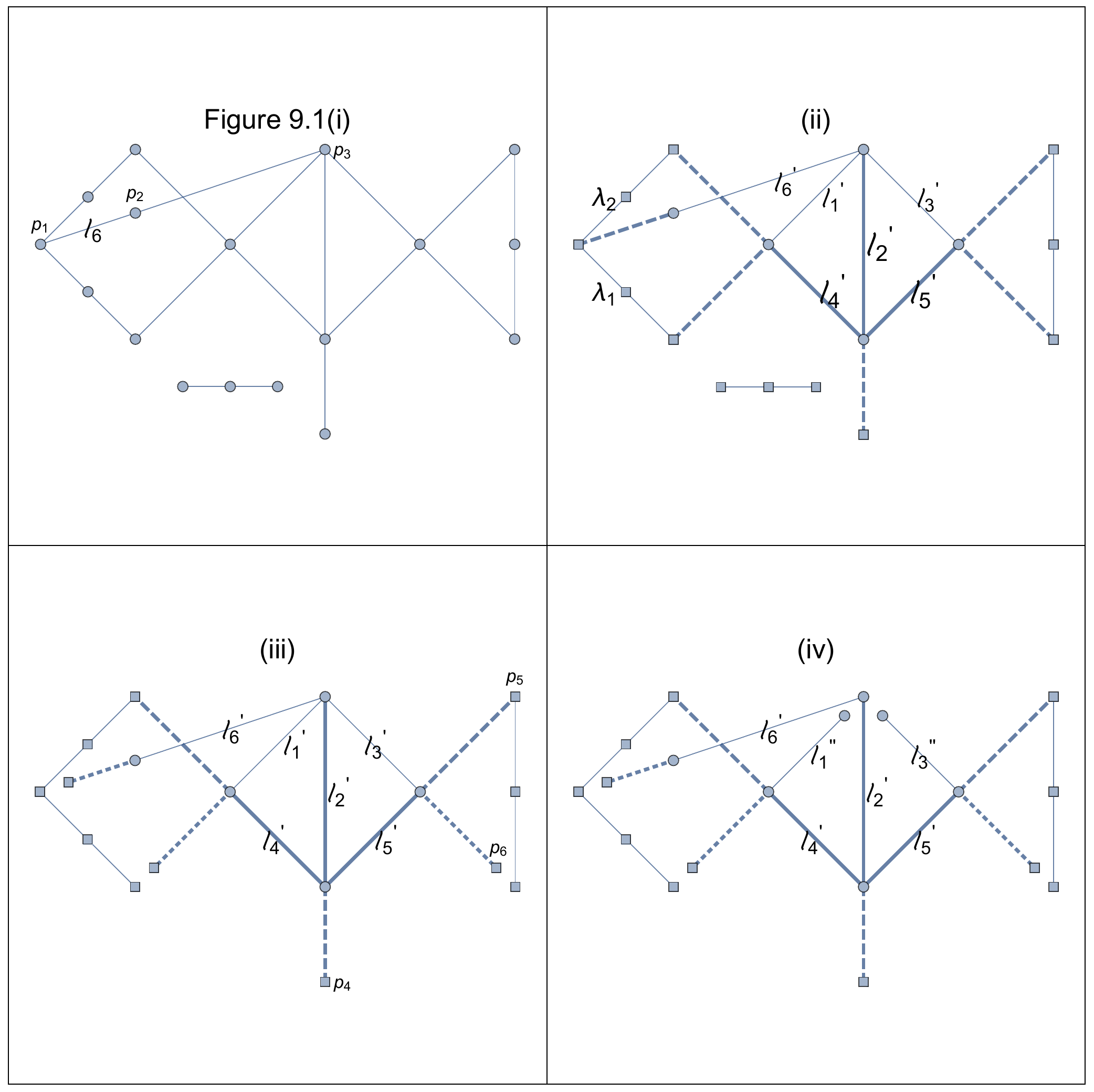} 
 \end{center}
 
 To illustrate, consider a hypothetical modular lattice $L$ with a (disconnected) base of lines $\B=(J(L),\LL)$ as shown in Figure 9.1(i).
 For later purposes one line $\ell_6=\{p_1,p_2,p_3\}$ of $\LL$ is showcased.
  Suppose $a\prec 1$ is such that $\B(a)$ decays in $s_a:=s([0,a])=4$ connected components; they have cardinalities 1,3,3,5 and their
   points are rendered by little squares in Figure (ii).
 
  What can be said about $\B(a,1)=(J(a,1),\LL(a,1))$? Since five $p\in J(L)$ are not squared, we have $|J(a,1)|=5$. Six lines of $\Lambda$ involve points of $J(a,1)$, call them $\ell_1$ to $\ell_6$. Their restrictions to $J(a,1)$ yield $\Lambda(a,1)=\{\ell_1',...,\ell_6'\}$ as shown in Fig.(ii).
  Each $\ell_i'$ in $\LL(a,1)$ 'points' (see dashed lines) to exactly one connected component of $\B(a)$. 
  Thus for instance $\ell_6'\cup\{p_1\}=\{p_2,p_3\}\cup\{p_1\}=\ell_6$.
  (However, not all connected components of $\B(a)$ are reached that way, i.e. not the horizontal 3-element line. Figures (iii) and (iv) are not relevant yet.)

 {\bf Theorem 9.1 [HW,5.2]:} {\it All localizations of bases of lines of  modular
 	lattices are connected.}
 
 \medskip
 {\it Proof.} Using induction on $\delta(L)$ we may consider a covering
 $a\prec 1$ and the localization $\B(a,1)$ of a given base of lines $\B=(J(L),\LL)$.
  If $a$ is the {\it only} coatom, then $\B(a,1)$ boils down to the singleton
 $\{1\}$ (which needs no lines to be connected). 
  Henceforth assume there are $m\ge 1$ 
 coatoms $b_1,\dots,b_m$ {\it apart} from $a$. Then $\B(a,1)=(J(a,1),\LL(a,1))$ with
 $J(a,1)=\bigcup_{i=1}^m J(ab_i,b_i)$. By induction, since $\delta([0,b_i])<\delta(L)$,
 all localizations $\B(ab_i,b_i)=(J(ab_i,b_i),\LL(ab_i,b_i))$
 are connected. Here the line sets $\LL(ab_i,b_i)=\{\ell\in \LL(a,1):\ \ol{\ell}\le b_i\}$ need not be disjoint, and their union is a proper subset of $\LL(a,1)$ iff  there is a line $\ell\in\LL(a,1)$ with $\ol{\ell}=1$. 
 
 So let us proceed to show that $\B(a,1)$ is connected.
 {\it First case:} $m=1$. Then $J(a,1)=J(ab_1,b_1)$, so $\B(a,1)=\B(ab_1,b_1)$
 is connected. {\it Second case:} $m\ge 2$.
 It suffices to show that, say, $\B(ab_1,b_1)$ is
 connected with $\B(ab_2,b_2)$. {\it Subcase (i):} $ab_1,ab_2,b_1 b_2$ are
 distinct (see Figure 9.2(A)).
 Putting $d:=ab_1 b_2$ one has $[d,1]\simeq B_8$, the 8-element Boolean lattice. Any $p\in J(d,b_1b_2)$
 lies in $J(ab_1,b_1)\cap J(ab_2,b_2)$, so $\B(ab_1,{b_1})$ and 
 $\B({ab_2},{b_2})$ have even a point in common. {\it Subcase (ii):}
 $ab_1=ab_2=b_1b_2=:d$ (see Figure 9.2(B)).
 Pick $p\in J(d,b_1)$ and $q\in J(d,b_2)$. Then
 $[p_*,p]$ and $[q_*,q]$ have $[a,1]$ as a common upper transpose. By
 Lemma 7.2 there is a $M_n$-element $x$ with $p+q=x$. Obviously $x\not\le a$, and so there is exactly one line $\ell'\in\LL(a,1)$ that belongs to the\footnote{A little extra thought (Ex.9A) shows that actually $[x_0,x]$ coincides with $[d,1]$.} line-interval $[x_0,x]$. Let $p',q'\in\ell'$ belong to the same atoms of $[x_0,x]$ as $p,q$ respectively. We conclude as follows that $p,q$ are connected in $\B(a,1)$: $p$ is connected to $p'$ in
 $\B(ab_1,b_1)$, line $\ell'$ brings us from $p'$ to $q'$, and $q'$ is connected to $q$ in $\B(ab_2,b_2)$.
 $\square$
 
 \begin{center}
 	\includegraphics[scale=0.8]{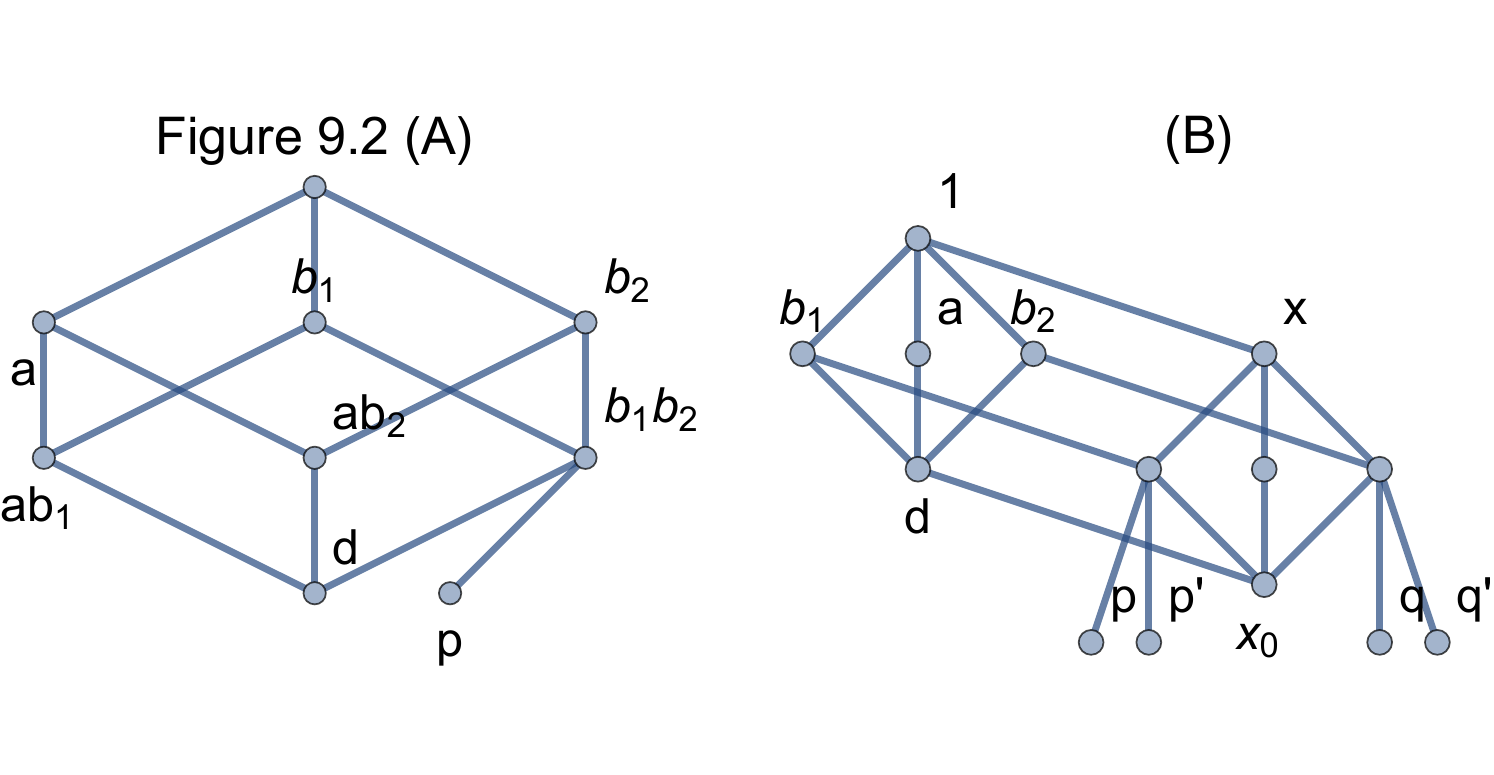} 
 \end{center}

{\bf 9.2}  Apart from  $j(L),\delta(L),s(L)$, here come some new  parameters. For each modular $L$ we define $i(L)$ as the number of line-{\it intervals} of $L$. Equivalently, $i(L)$ is the number of $M_n$-elements (=line-tops), or the number of lines in any base of lines. For instance $i(L_1)=3$ (Figure 3.1). Further $o(L)$ is the maximum value of $n-1$ achieved by the $M_n$-elements, so $o(L_1)=3$. Using $n-1$ rather than $n$ will enhance the look of upcoming formulas.
Although we defined $M_n$-elements only for $n\ge 3$, it is convenient here to set $o(L):=1$ if $L$ is distributive.
 Finally $\mu(L)$ is the sum of all $n$'s that occur in the $M_n$-elements of $L$. 

{\bf 9.3} After these numerical parameters, let us introduce two important kinds of lattices.
First, a modular lattice $L$ is \underbar{acyclic} if all its bases of lines are acyclic, i.e. without cycles in the sense of PLSes. Likewise $L$ is \underbar{cyclic} if it is not acyclic, i.e. at least one BOL has cycles.
 
 {\bf Theorem 9.2: }
{\it For each finite modular lattice $L$ the following holds.
	\begin{enumerate} 
	\item[(a)]  $j(L)\le \mu(L)-i(L)+s(L)$.
    \item[(b)] $L$ is acyclic $\LRa\ j(L)= \mu(L)-i(L)+s(L)$.
\end{enumerate} } 

The proof is left as Exercise 9B.  It follows at once from Theorem 9.2 (b) that 
in a {\it finite} modular lattice $L$ the acyclicity of one BOL is equivalent to the acyclicity of all BOL's, i.e. to the acyclicity of $L$. In contrast, for {\it any} infinite (but FL) modular lattice the identity in (b) boils down to the void statement $\infty=\infty$. Nevertheless, the author speculates that the answer to this question is "yes":

 {\bf Open Question 1: } 
{\it  Is it true that also for infinite (but FL) modular lattices the acyclicity of one BOL is equivalent to the acyclicity of all BOL's? }

{\bf 9.3} Motivated by Figure 5.2 (that illustrates the Triangle Axiom of projective spaces) we say that a BOL $\B=(J(L),\Lambda)$ of any modular lattice $L$ contains a \underbar{triangle configuration} if there are distinct $\ell_1,...,\ell_4\in\Lambda$ such that $\{\ell_1,\ell_2,\ell_3\}$ form a triangle (say with points of intersection $s,q,r$) and $\es\neq\ell_4\cap\ell_i\not\s \{s,q,r\}$ for $i=1,2,3$.

{\bf Theorem 9.3: }{\it  Suppose the base of lines $\B$ of the modular lattice $L$ contains a triangle configuration. Then 
some localization $\B(a,b)$ has a cycle.}

{\it Proof.} Our argument elaborates the direction $(v)\Ra (iv)$ in [HPR,Thm.7.2]. Adhering to the notation of Figure 5.2 put $\ell_1=\langle s,p_1\rangle,\ell_2=\langle s,p_2\rangle,\ell_3=\langle p_1,p_2\rangle$. The main idea is to show that for some suitable covering $(a,b)$ the cycle $(\ell_1,\ell_2,\ell_3)$ induces a cycle $(\ell_1',\ell_2',\ell_3')$ in $\B(a,b)$. Claim (1) will pave the way to define $(a,b)$.

\begin{enumerate}
	\item [(1)]{\it The $M_n$-element $u:=q+r=q+p_3=r+p_3$ is such that $s\not\le u$.  }
\end{enumerate}

By way of contradiction, assume that $s\le u$. Hence, $s$ is below one of the co-atoms of $u$, say $s\le u_1$. (Note that $s\neq u=q+r$ by join-irreducibility.)  {\it Case 1:} $s\not\le u_0$. Since $u_1$ coincides with at most one of $u_0+q,u_0+r,u_0+p_3$ (these are coatoms of $u$ by Lemma 7.1), we conclude wlog that $u_1, u_0+q, u_0+r$ are distinct. But this implies $q+s=q+r=u$, which contradicts $\langle q,s\rangle\neq\langle q,r\rangle$. {\it Case 2:} $s\le u_0$. Then $p_1< q+s\le q+u_0$. Assuming $p_1\not\le u_0$ would imply $p_1+p_3=u$, which contradicts $\ol{\ell_3}\neq u$.
Therefore $p_1\le u_0$. Yet this together with $s\le u_0$ leads to the contradiction $q<p_1+s\le u_0$. This finishes the proof of (1).

According to (1) we have $s\not\le u$, which implies $su\le s_*\prec s$. From $[su,s]\simeq [u,u+s]$ follows that
$a:=u+s_*\prec u+s=:b$. It remains to show that (i) $s,p_1,p_2\le b$, but (ii) $s,p_1,p_2\not\le a$.
As to (i), from $s< r+p_2 \le u+p_2$ and $p_2<r+s\le u+s$ follows $u+p_2=u+s=b$. Similarly from $s<q+p_1\le u+p_1$
and $p_1<q+s\le u+s$ follows $u+p_1=u+s=b$. As to (ii), $a=u+s_*\prec u+s$ forces $s\not\le a$. Assuming $p_1\le a$ leads to the contradiction $s<p_1+q\le a+u=a$; and assuming $p_2\le a$ leads to  $s<p_2+r\le a+u=a$. $\square$

{\bf 9.4} A modular lattice $L$ is \underbar{locally acyclic} if for all bases of lines $\B$  all localizations $\B(u,v)$ are acyclic. It is evident that 'acyclic $\Ra$ locally acyclic'. Although locally acyclic lattices can be cyclic (Ex.9C), by Theorem 9.3 they cannot contain triangle configurations.
We now argue that even in the 'nasty' case\footnote{The characterization of 'nasty' will be postponed to Part C.} that
  a localized BOL is cyclic, there always remains a glimmer of acyclicity:

\begin{enumerate}
	\item [(2)] {\it Fix any c-simple modular lattice $L$. Hence by Thm.7.3
		 each base of lines $\B$ is connected. Further fix any covering $(a,1)$ and let $\B(a,1)$ be its localization. Suppose $\B(a)$ has $s_a:=s([0,a])$ connected components $\B_i$. Since $\B$ is connected, we can pick for each $\B_i$ any line that connects $\B_i$ with $\B(a,1)$.
		  This yields a $s_a$-element set  $T\s\LL(a,1)$, and so $|\LL(a,1)|\ge s_a$. Since moreover $T$  can be proven\footnote{ The lengthy argument (the knowledgeable reader is invited to shorten it) is given in [Wild, Algebra Universalis 35 (1996),p.117]. Other than mistakenly stated in [HW1,p.31,line 5] there is {\it no obvious}  reason why $|\LL(a,1)|\ge s_a$ should imply
			$j(a,1)\ge s_a+1$. In fact, e.g.  $|\LL(a,1)|=6$ but $j(a,1)=4$  for $L=L(\Z_2^3)$; see Ex.9I.} to be {\bf acyclic}, (4) in Section 5 implies that $|\bigcup T|\ge s_a+1$. A fortiori $j(a,1)\ge s_a+1$.}
\end{enumerate}

To illustrate, imagine that $\B$ in Figure 9.1(A),(B) was connected by mentally erasing the isolated 3-element line. Instead of $s_a=4$ we now have $s_a=3$. One of four ways to choose  $s_a$ lines is $T:=\{\ell_2',\ell_4',\ell_5'\}$, which is rendered boldface in (B).

{\bf 9.5}
Recall from Sec. 5.4 that $r^*(\B)$ is the number of point-splittings required to turn the partial linear space $\B$ into an acyclic PLS with the same number of connected components. 
 
\bigskip

 {\bf Theorem 9.4 (adding point-splittings to [HW,6.4]): }
 {\it Let $L$ be a modular lattice with a base of
 lines $\B=(J(L),\LL)$.
 \begin{enumerate}
 	\item[(a)] Generally $i(L)\ge\delta(L)-s(L)$ and $j(L)\ge 2\delta(L)-s(L)$.
 	\item[(b)]If moreover $o(L)\le 2$,
 	then 	$j(L)\ \ge\  2i(L)+s(L)-r^*(\B)\ \ge\ 2\delta(L)-s(L)$.
 	\item[(c)] Suppose $L$ is LOCALLY ACYCLIC.\\ Then $i(L)=\delta(L)-s(L)+r^*(\B)$
 	and	$j(L)\ge i(L)+\delta(L)$.\\ If moreover $o(L)\le 2$, then $j(L)=i(L)+\delta(L)$.
 	\item[(d)] Suppose $L$ is ACYCLIC.\\ Then $i(L)=\delta(L)-s(L)$.\\
 	 If	moreover  $o(L)\le 2$, then
 	$j(L)=2\delta(L)-s(L)$.  
 \end{enumerate}}

 \medskip\rm
 {\it Proof.} {\bf (a).} We use induction on $\delta(L)$. 
 
 {\it First case:} $s:=s(L)\ge 2$. Let
 $L_1=L/\theta_1,\dots,L_s=L/\theta_s$ be the c-simple factors with corresponding bases of lines
 $\B_j=(J_j,\LL_j)$. As to the first inequality, since $\delta(L_j)<\delta(L)$ one has
 $i(L_j)\ge\delta(L_j)-1\ (1\le j\le s).$
  Together with 
  $\delta(L)=\sum_{j=1}^s\delta(L_j)$ (Theorem 3.3) follows
 $i(L)=\sum_{j=1}^s i(L_j)\ge \delta(L)-s(L)$.
As to the second inequality in (a), again by induction $j(L_i)\ge 2\delta(L_i)-1$. By (10) in Section 3
we have $j(L)=\sum_{i=1}^s j(L_i)$. It follows that
  $j(L)\ge 2\delta(L)-s(L)$.
 \smallskip

 {\it Second case:} $s=1$. Consider a coatom $a\prec 1$ and the induced
 base of lines $\B(a)=(J(a),\LL(a))$. According to Theorem 
 7.3 it has $s_a:=s([0,a])$
 connected components $\B_i'$. Since $\delta(a)<\delta(L)$, induction (anchored for $\delta(L)=1$, check) gives
 
 $ (3)\quad i([0,a])\ge \delta(a)-s_a=\delta(L)-1-s_a,\ as\ well\ as\ j(a)\ge 2\delta(a)-s_a=2\delta(L)-2-s_a$. 
 
 As to the claim $i(L)\ge \delta(L)-1$ of part (a), by (2) we have $|\LL(a,1)|\ge s_a$ which, jointly with the first part of (3), yields
 $i(L)=i([0,a])+|\LL(a,1)|\ge i([0,a])+s_a   \ge\delta(L)-1$. As to the claim $j(L)\ge 2\delta(L)-1$, by  fact (2) we have $j(a,1)\ge s_a+1$. This together with the second part of (3) gives
 $j(L)=j(a)+j(a,1)\ge   j(a)+s_a+1\ge 2\delta(L)-1$.

{\bf (b).}  The first claim $j(L)\ge 2i(L)+s(L)-r^*(\B)$ 
 is equivalent to $j(L)+r^*(\B)\ge 2i(L)+s(L)$.
 Applying $r^*(\B)$ point-splittings to $\B$ yields an acyclic
  partial linear space $\B^0$ with $j(L)+r^*(\B)$ points and $c(\B^0)=c(\B)=s(L)$.
  Being acyclic each connected component $\B^k=(J^k,\LL^k)$ of $\B^0$ satisfies $|J^k|\ge 2|\LL^k|+1$.
  In view of $|\LL^k|=i(L_k)$ and Thm.7.3 summing over $k=1,..,s(L)$ yields $j(L)+r^*(\B)\ge 2i(L)+s(L)$.

 The second claim $2i(L)+s(L)-r^*(\B)\ge 2\delta(L)-s(L)$ in (b) is
 equivalent to the claim

 $(4)\quad r^*(\B)\le 2i(L)+2s(L)-2\delta(L)$.
 
  Establishing the latter will be more subtle and, different from above, will exploit the assumption $o(L)\le 2$. We will use again the hypothetical BOL $\B$ in Fig.9.1 (with the mentally erased isolated line) to fix ideas.

 {\it First Case: }$s=s(L)\ge 2$. Since $\delta(L_k)<\delta(L)$, induction gives $r^*(\B^k)\le 2i(L_k)+2-2\delta(L_k)$. Upon summing up (and invoking Theorems 3.3 and 7.3) we obtain
 $r^*(\B)\le 2i(L)+2s(L)-2\delta(L)$. 
 
 {\it Second Case:} $s=1$. Picking $a\prec 1$ claim (4) becomes
 
 $(5)\quad r^*(\B)\le 2i(L)+2-2\delta(L)=2i(L)-2\delta(a).$
 
 Since by induction (4) holds for $[0,a]$, we have
 
 $(6)\quad r^*(\B(a))\le 2\vv\LL(a)\vv+2s_a-2\delta(a).$

 Subtracting (6) from (5) it remains  to show that
 
 $(7)\quad r^*(\B)- r^*(\B(a))\le 2\vv\LL(a,1)\vv-2s_a.$

For notational convenience we henceforth apply $r^*$ directly to line-sets rather than PLSes, say $r^*(\LL(a)$ instead of $r^*(\B(a))$. This is justified by the fact hat the cardinality of an acyclifier of a PLS does not depend on the number of isolated points of the PLS.
 Thus, since inequality (7) amounts to $r^*(\LL)\le r^*(\LL(a))+(|\LL(a,1)|-s_a)+(|\LL(a,1)|-s_a),$
  it suffices to show that $r^*(\LL)$ can be written as
 
$(8)\quad r^*(\LL)=r^*(\LL(a))+r^*(\LL(iii))+\beta$

in such a way that

$(9)\quad  \beta= |\LL(a,1)|-s_a,$
as well as

$(10)\quad  r^*(\LL(iii))\le |\LL(a,1)|-s_a.$

 As to (8), we define $\beta$ as the number of point-splittings necessary to destroy {\it mixed} $\LL$-cycles, i.e. having junctions $\beta$oth in $J(a)$ and $J(a,1)$. (Examples of mixed cycles in Figure (ii)  are $\{\ell_6,\lambda_2,\ell_4,\ell_2\}$ and $\{\lambda_2,\lambda_1,\ell_1,\ell_4\}$.) Generally these $\beta$ point-splittings can be performed by detaching $\beta$ lines from $\LL\setminus\LL(a)$ from points in $J(a)$; let  $\LL(iii)$ be the arising  line-set. (As can be seen in Fig.(iii), in our example $\beta=3$.)   We next invest $r^*(\LL(iii))$ point-splittings for destroying all cycles only involving lines from $\LL(iii)$. (In Fig.(iii) such a cycle e.g. is $\{\ell'_2\cup\{p_4\},\ell'_5\cup\{p_5\},\ell'_3\cup\{p_6\}\}$.) Finally $r^*(\LL(a))$ point-splittings are invested to destroy all $\LL(a)$-cycles. After that $\LL$ is completely acyclified. (In our example $r^*(\LL(iii))=2$ and $r^*(\LL(a))=0$, and Fig. (iv) shows the completely acyclified line-set triggered by $\LL$.)
 
 %akin to 5.4.1 we choose an acyclifier of $\LL$ that suits us.  For starters we invest $r^*(\LL(a))$ point-splittings to  destroy all $\LL$-cyles that  only involve lines from $\LL(a)$. (As seen in Fig.9.1(B), in our toy example $\LL(a)$ is acyclic, and so $r^*(\LL(a))=0$.) Next  invest $r^*(\LL\setminus\LL(a))$ point-splittings to destroy all $\LL$-cyles  that only involve lines from $\LL\setminus\LL(a)$. Although $\LL\setminus\LL(a)$ is not quite the same as $\LL(a,1)$, it is clear that $r^*(\LL\setminus\LL(a))=r^*(\LL(a,1))$. (From Fig.(C) one sees that in our example $r^*(\LL(a,1))=2$.) Finally we define $\beta$ as the number of point-splittings required to destroy all $\LL$-cyles that involve  lines from {\it both} $\LL(a)$ and $\LL(a,1)$. (Figure (D) shows that in our example $\beta=3$. For instance the three shown point-splittings destroy the "mixed" cycles $\{\ell_6,\lambda_1,\ell_4,\lambda_2\}$ and $\{\lambda_1,\lambda_2,\ell_1,\ell_4\}$. (If we were pedantic, we wouldn't speak of $\ell_1$ here, but rather of what has become of $\ell_1$ after a previous point-splitting.))   
 
 {\it Proof of (9).}  Let $\ol{T}\s\LL$ be the set of extensions of lines from $T$. (Thus $\ol{T}=\{\ell_2,\ell_4,\ell_5\}$ in Fig.(ii).) Further put $\ol{\LL}=(\LL\setminus\LL(a))\setminus\ol{T}$. (Thus $\ol{\LL}=\{\ell_1,\ell_3,\ell_6\}$ in Figure (ii).) We claim that each $\ell\in\ol{\LL}$ must undergo exactly one point-splitting since otherwise mixed cycles remain. (Once this is established, we conclude $\beta=|\ol{\LL}|=|\LL\setminus\LL(a)|-|\ol{T}|=|\LL(a,1)|-s_a$ as claimed.)
  Indeed, each $\ell\in\ol{\LL}$ connects $\B(a,1)$ with some connected component $K$ of $\B(a)$. Since there is a line $\ell_0\neq\ell$ in $\ol{T}$ that leads back from $K$ to $\B(a,1)$, and since $\B(a,1)$ itself is connected (Thm.9.1), there is a cycle $C$ which has junctions both in $J(a)$ and $J(a,1)$. (Notice that $C$ may  lack {\it lines} from either $\LL(a)$ or  $\LL\setminus\LL(a)$, but of course not from both.) Let $\LL(iii)$ be the line-set obtained from $\LL\setminus\LL(a)$ upon "detaching" each $\ell\in\ol{\LL}$ from its unique point in $J(a)$. Of course $|\LL(iii)|=|\LL\setminus\LL(a)|$. Suppose there still was a 'mixed' cycle $C$ using lines from $\LL(iii)\cup\LL(a)$. (The meaning of 'mixed' has slightly changed now: $C$ must have junctions in both $J(a)$ and the remaining point-set $J'$ of cardinality $j(a,1)+\beta$.)
   Then $C$ connects $J'$ to some component $K$ of $\B(a)$, and this necessarily by virtue of a line $\ell_1\in \ol{T}$ (because other connecting lines have been detached). In order for $C$ to come back from $K$ to $J'$, it must  use {\it another} line $\ell_2\in\ol{T}$. This is impossible since by construction $\ol{T}$ contains exactly one line connecting $K$ and $J'$. Thus the line-set $\LL(iii)\cup\LL(a)$ obtained so far has no more mixed cycles. Hence each remaining cycle either has all its junctions in $J'$ (and whence all its lines in $\LL(iii)$) or it has all its junctions in $J(a)$ (and whence all its lines in $\LL(a)$).

 {\it Proof of (10)} If $\LL(iii)$ contains a cycle $C$, apply one point-splitting to one of its lines $\ell$. This destroys $C$ without increasing the number of connected components.
 
  Furthermore, since $o(L)\le 2$ entails $|\ell|=2$, {\it all} cycles that utilized $\ell$ are now destroyed.
  
  If the new PLS still has cycles $C_0$, continue in this manner until the PLS is acyclic. In order to get an upper bound for $r^*(\LL(iii))$ we choose particular point-splittings. Whenever a cycle $C_0=(\ell_1,\ldots,\ell_n)$ needs to be destroyed, we pick a line $\ell_i$ that doesn't belong to $\ol{T}$. This can be done because $\{\ell_1,\ldots,\ell_n\}\s \ol{T}$ is prevented by the acyclicity of $T$ (see (2)). From this and from "Furthermore..." above follows that\\ $ r^*(\LL(iii))\le |\LL(iii)|-s_a=|\LL\setminus\LL(a)|-s_a= |\LL(a,1)|-s_a$.
 
% {\it Proof of (10).}  Let $\LL_0$ be the set of all extended lines $\ell=\ell'\cup\{p_\ell\}$, where $\ell'$ ranges over $\LL(a,1)\setminus T$. (In our example $\LL_o=\{\ell_1,\ell_3,\ell_6\}$ where $\ell_1',\ell_3',\ell_6'$ are from Figure (C).) We claim that each $\ell\in\LL_0$ must undergo a point-splitting since otherwise cycles $C$ remain. Indeed, $\ell$ connects $\B(a,1)$ with some connected component $K$ of $\B(a)$. Since there is a line $\ell_0\in T$ that leads back from $K$ to $\B(a,1)$, and since $\B(a,1)$ itself is connected (Thm.9.1), there must be cycles $C$. Let $\ol{\B}$ be the PLS obtained after detaching all $\beta$ many lines of $\LL_0$ (see Fig.(D)). If there still was a cycle $\ol{C}$ in $\ol{\B} $, then $\ol{C}$ connects $\B(a,1)$ to some component $K$ of $\B(a)$, and this necessarily by virtue of a line $\ell_1\in T$ (because no other ones are left). In order for $\ol{C}$ to come back from $K$ to $\B(a,1)$, it must use {\it another} line $\ell_2\in T$. This is impossible since by construction $T$ contains exactly one line connecting $K$ and $\B(a,1)$. It follows that $r^*(\ol{B})=\beta=|\LL(a,1)|-s_a$.  

 {\bf (c).} The first claim is equivalent to claiming $r^*(\B)=i(L)+s(L)-\delta(L)$. We induct on $\delta(L)$. {\it First case:} $s(L)\ge 2$. This is easily handled akin to similar situations in (a) and (b).

 {\it Second case:}
 $s(L)=1$. Consider a co-atom $a\prec 1$ and the corresponding
 PLSes $\B(a),\
 \B(a,1)$. With $L$ also $[0,a]$ is locally acyclic, and so by induction
 $r^*(\B(a))=\vv\LL(a)\vv+s_a-\delta(a)$, where again $s_a:=s([0,a])$.
 Recall from part (b) that generally $r^*(\B)=r^*(\B(a))+r^*(\B(a,1))+\beta$, where $\beta=|\LL(a,1)|-s_a$.
 Because $r^*(\B(a,1))=0$ by the very definition of local acyclicity, we conclude that
 
  $r^*(\B)=(\vv\LL(a)\vv+s_a-\delta(a))+0+(\vv\LL(a,1)\vv-s_a)=
 i(L)+1-\delta(L),$
 
 which was to be proven.
 As to the second claim in (c), using the first statement of (c) in the form $i(L)-r^*(\B)=\delta(L)-s(L)$, as well as (b), yields
 
 $j(L)\ge 2i(L)+s(L)-r^*(\B)=i(L)+(i(L)-r^*(\B))+s(L)=\\ i(L)+(\delta(L)-s(L))+s(L)=i(L)+\delta(L).$
 
 Assume that moreover $o(L)\le 2$, so all lines in $\LL(a,1)$ have cardinality 2. Since $\B(a,1)$ is connected by Theorem 9.1 and acyclic by assumption, we have $j(a,1)=|\LL(a,1)|+1$.
 Furthermore by induction 
 $j(a)=\vv\LL(a)\vv+\delta(a)$. Therefore $j(L)=j(a)+j(a,1)=|\LL|+1+\delta(a)=i(L)+\delta(L)$.
 
 {\bf (d).} By assumption $r^*(\B)=0$, and so $i(L)=\delta(L)-s(L)$ by part (c). If moreover $o(L)\le 2$ then 
 $j(L)=i(L)+\delta(L)$ again by part (c). Hence $j(L)=i(L)+\delta(L)=(\delta(L)-s(L))+\delta(L)=2\delta(L)-s(L)$.
 $\square$

 \bigskip

 By part (c) in each locally acyclic lattice 
 $r^*(\B)=i(L)+s(L)-\delta(L)$ is independent of the actual
 base of lines $\B$. 
  In contrast, part (b) of Theorem 9.4 only gives  upper and a lower bounds for the range of  
 $r^*(\B)$, namely
 
 $(11)\quad 2i(L)+s(L)-j(L)\le r^*(\B)\le 2i(L)+2s(L)-2\delta(L)$
 
 Notice that the upper bound is exactly twice the value of $r^*(\B)$ in the locally acyclic case.
 
  {\bf Open Question 2:} {\it Is $r^*(\B)$ the same for all BOL's $\B$ of any fixed modular lattice?}
  
  \vspace{1.2cm}

{\bf 9.6}
Apart from cycles in BOL's we introduce another kind of cycle. For $M_n$-elements $x,y$ we say\footnote{This is a handy  definition for a relation that apparently was so far nameless.} that $x$ is \underbar{ssmaller} than $y$, written $x<^* y$, if $x<y$ but $x\not< y_0$. Likewise 
$x$ is \underbar{llarger} than $y$, written $x>^* y$, if $y<^* x$. Finally $x,y$ are \underbar{ccomparable}
if one of the two takes place. Similar to [Mitschke-Wille, \"Uber freie modulare Verb\"ande $FM(M_D^q)$, Contributions to General Algebra, Klagenfurt 1978] we call a sequence $(x,y,...,z)$ of $M_n$-elements \underbar{cycle of $M_n$-elements} if all successive elements (including $z,x$) are ccomparable.

\begin{center}
	\includegraphics[scale=0.8]{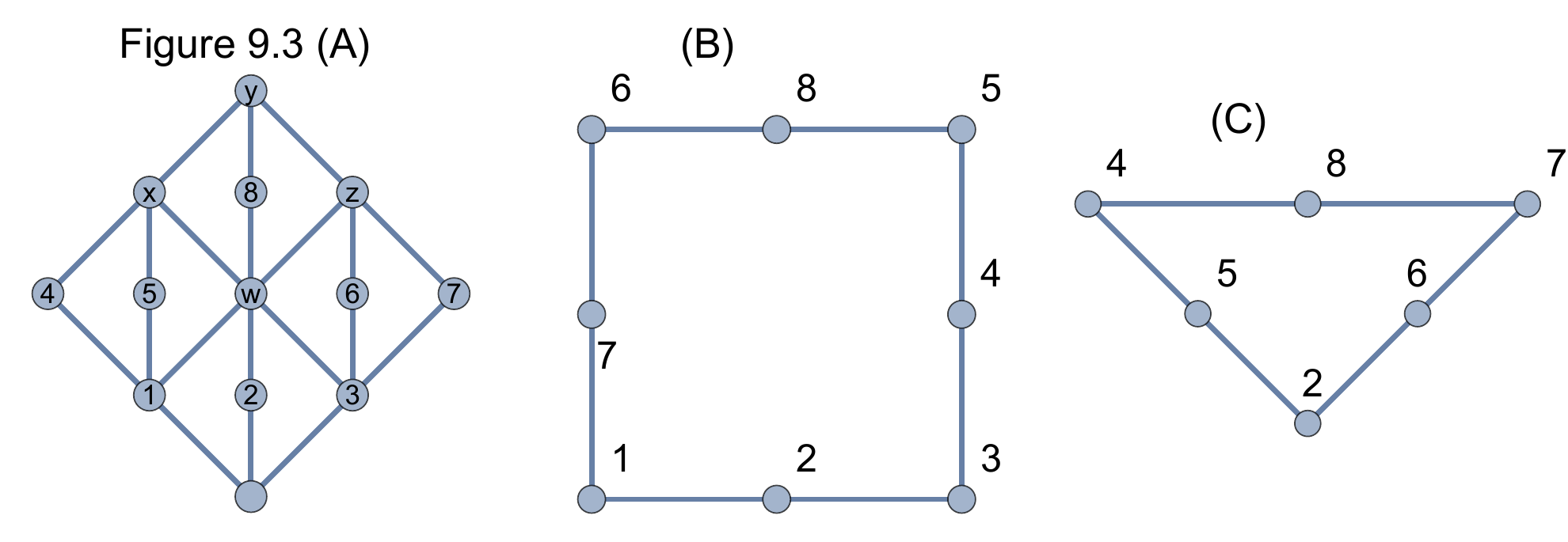} 
\end{center}

The lattice  in Figure 9.3 (A) has the cycle of $M_n$-elements $(w,x,y,z)$. The cycle of four lines in (B) is such that the corresponding line-tops are the four $M_n$-elements. Unfortunately the relation between the two kinds of cycles isn't always that straightforward, as witnessed (see (C)) by the cycle of three lines, whose corresponding line-tops $(x,y,z)$ do {\it not} constitute a cycle of $M_n$-elements.
Let us shed more light on the matter. 

{\bf 9.6.1} By definition of $M_n$-element, when $x<^* y$ then\footnote{Notice that equality $x=y_i$ e.g. happens in Fig.9.3(A). } $x\le y_i$ for exactly one  $y_i\prec y$. That's because $x\le y_i,\!y_j$ gives the contradiction $x\le y_i y_j=y_0$.  Clearly $x+y_0=y_i$, and so $xy_0\prec x$, say $xy_0=x_k$. In other words, $(x_k,x)\nearrow (y_0,y_i)$. Conversely, if such a relation exists between an \underbar{upper covering} of $[x_0,x]$ and a \underbar{lower covering} of $[y_0,y]$ (obvious definitions), then $x<^* y$. So much about {\it two} subsequent elements in a $M_n$-cycle.

\begin{center}
	\includegraphics[scale=0.9]{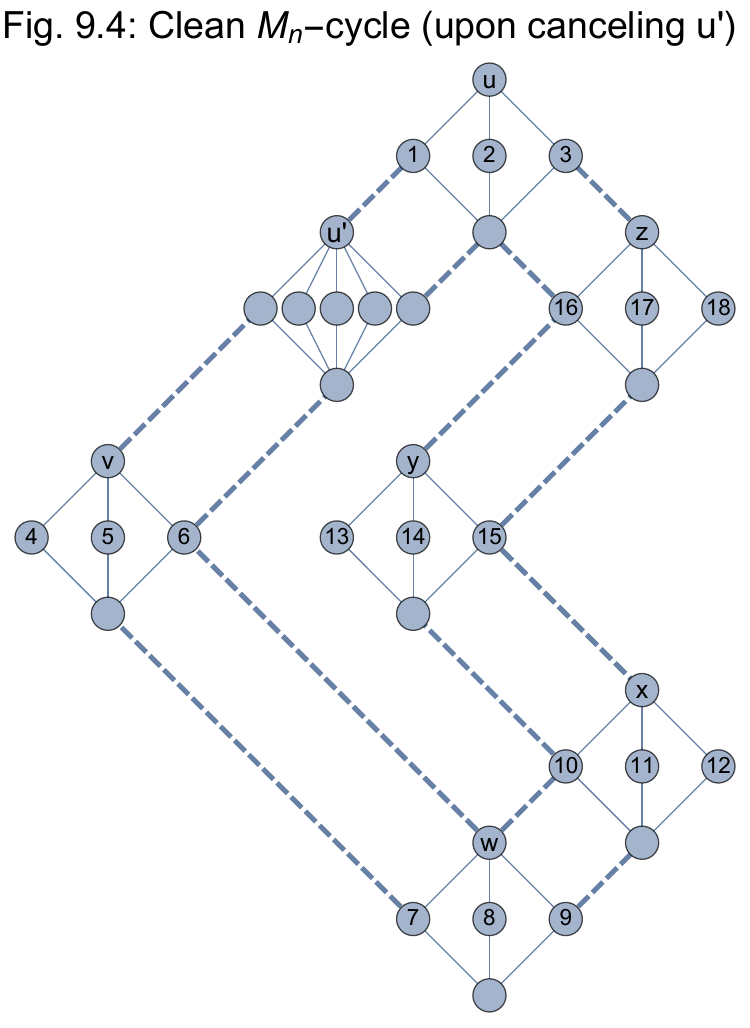} 
\end{center}

This cycle of seven $M_n$-elements (without surrounding lattice) will be large enough to detect general techniques to come up with an induced cycle of lines.
According to Figure 9.4 we have

$(12)\quad u>^*u'>^* v>^*w<^*x<^*y<z^*<^*u.$

Apart from $u'$ all $M_n$-elements have $n=3$. This makes it easier for the eye. Furthermore, the minimum value $n=3$ is big enough for showing how the junctions of the lines in our upcoming cycle of lines must be chosen.
 As to $u'$, whenever two steps in the same direction, such as $u>^*u'>^* v$, can be\footnote{Given $u>v$, it holds that $u>^*v$ iff $u_0\not\ge v$. This takes place here, but e.g. in $w<^*x<^*y$ we have $w<y_0$, and so $w\not<^* y$.  } \underbar{condensed} to $u>^*v$, then dropping the middle element $u'$ still leaves a $M_n$-cycle. After dropping all such $M_n$-elements, the shorter $M_n$-cycle has a zig-zag pattern in the sense that none of the four "directions" NE,NW,SE,SW gets repeated immediately.
 
  Let us construct a cycle of lines whose line-tops yield the cleaned $M_n$-cycle of Fig.9.4. Any line $\ell_1$ with $\ol{\ell_1}=u$ has the form $\{1',2',3'\}$
where $1'\in J(u_0,1),2'\in J(u_0,2),3'\in J(u_0,3)$. Similar notation will be used for the other $M_n$-elements in Fig.9.4.
Can $\ell_2$ with $\ol{\ell_2}=v$ be chosen in such a way that $\ell_1\cap\ell_2=\{p\}$ with $p\in J(u_0,2)$?
Then $p\le 2\wedge v=6\le u_0$, yet no element of $\ell_1$ is $\le u_0$. Similarly $p\in J(u_0,3)$ is impossible.
We see that in order for $\ell_2$ to connect with $\ell_1$ via $p$ we must either have $p\in J(v_0,4)\s J(u_0,1)$ or $p\in J(v_0,5)\s J(u_0,1)$. Choosing (say) $p=4'\in J(v_0,4)$ yields
 $\ell_1=\{\ul{4'},2',3'\}$. We will always underline the element that provides the junction to the next line. Similar reasoning yields $\ell_2=\{4',5',\ul{8'}\}$ as one possibility. Continuing in this manner we e.g. get
$\ell_3=\{\ul{7'},8',9'\}, \ell_4=\{7',11',\ul{12'}\}, \ell_5=\{12',\ul{13'},14'\}, \ell_6=\{13',\ul{17'},18'\}.$
Now we have literally come full circle, i.e. the next line $\ell_1$ has been fixed {\it already} as $\ell_1=\{4',2',3'\}$. Fortunately we can simply redefine the old $\ell_1$ as $\ell_1:=\{4',2',17'\}$ and thus get a cycle of lines.

But what  to do if instead of $3'$ we were somehow forced to substitute $4'$ in $\ell_1$? This cannot be done without breaking the link between $\ell_1$ and $\ell_2$, and this may have further repercussions. Perhaps that alley will be persued in a future version of Part B. 

%\begin{center}
%	\includegraphics[scale=0.9]{ModularFig9p5} 
%\end{center}

{\bf 9.6.2} For the time being we simply focuse   on \underbar{clean} $M_n$-cycles, which by definition are as follows:

\begin{enumerate}
	\item [(i)] Whenever three succesive $M_n$-elements $v,u,z$ stand in the relation $v<^*u>^*z$, it is forbidden that $v,u,z$ are {\it mutually} comparable, and that $(v_i,v)\nearrow (u_0,u_j)\nwarrow (z_k,z)$.
	\item [(ii)] Whenever three succesive $M_n$-elements $v,u,z$ stand in the relation $v>^*u<^*z$, it is forbidden that $v,u,z$ are {\it mutually} comparable, and that $(v_0,v_i)\nearrow (u_j,u)\nwarrow (z_0,z_k)$.
		\item [(iii)] No elements in the $M_n$-cycle get repeated.
\end{enumerate}

(We allow that parts of a clean cycle can be condensed.) The following  is clear from the deliberations in 9.6.1:

{\bf Theorem 6.5. }{\it If $L$ is a modular lattice which has a clean $M_n$-cycle, then there is at least one base of lines  with a cycle.   }

If $L$ is finite, then by Thm. 9.2 a clean $M_n$-cycle yields cycles in {\it all} BOL's of $L$.

{\bf Open Question 3:} {\it Does each $M_n$-cycle induce a clean $M_n$-cycle? If not, does it at least induce a $\B$-cycle in an appropriate BOL $\B$?}
\bigskip

{\bf 9.7} Let $a\prec 1$ be a co-atom in a modular lattice and let $C=(x,y,...,z)$ be a $M_n$-cycle which has
all its $M_n$-elements $\not\le a$. (Akin to localizations  $\B(a,1)$ in 9.1 each $M_n$-element has a unique coatom which is $\le a$.)  Further assume that  no three successive $M_n$-elements are mutually comparable. We postpone to a later version of Part B the argument  that the shape of $C$ must be as in Figure 9.5 below. Each circle labelled '$\le a$' indicates the location for the unique $M_n$-coatom  that is $\le a$.

%Argument is in cold coffee

\begin{center}
	\includegraphics[scale=1.2]{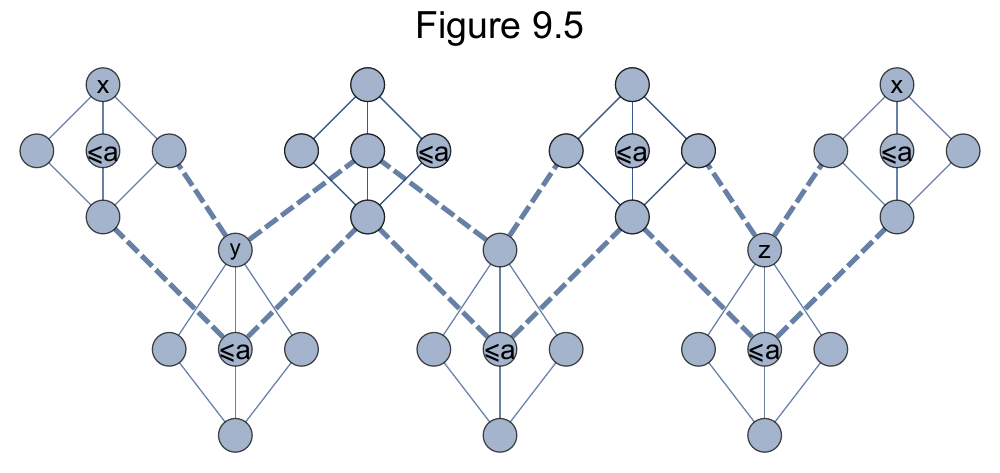} 
\end{center}

Any pair of parallel dashed lines can collapse. For instance for each co-atom $a$ of $L(\Z_2^3)$ {\it every} such pair collapses.

 \centerline{\bf 9.8 Exercises}

 {\bf Exercise 9A.} Show that the two intervals $[x_0,x]$ and $[d,1]$ occuring in the proof of Theorem 9.1 coincide.
 
 {\it Solution to Ex.9A.} If we can show that $[x_0,x]$ transposes up to $[d,1]$, then $[d,1]$ cannot be a line-interval, unless $[d,1]=[x_0,x]$. So suppose one had $d+x\le z$ with $d\prec z\prec 1$. We can assume that $z\neq b_1$ (otherwise $z\neq b_2$). But this yields the contradiction $p\le b_1 x\le b_1 z=d$. Therefore $d+x=1$.
 From $p_*,q_*\le d$ follows $x_0=p_*+q_*\le d$, and so $x_0\le dx$. From $\delta(dx)=\delta(d)+\delta(x)-\delta(d+x)=\delta(x)-2$ follows $x_0= dx$.
 
 {\bf Exercise 9B.} Prove Theorem 9.2.
 
 {\it Solution to Ex.9B.}  As to (a) in Thm.9.2, assume first that $L$ is c-simple. If $t=i(L)$ and the $t$ line-tops of $L$ (in any order) are a $M_{n_1}$-element, a $M_{n_2}$-element, ..., a $M_{n_t}$-element, then $\mu(L):=n_1+n_2+\cdots+n_t$. Let $\ell_1,...,\ell_t$ be any lines having these line-tops. By (2) in Sec. 5 we can assume that there are (not necessarily distinct) join irreducibles $p_i\in\ell_i\cap\ell_{i+1}$. Inducting on $i$ it is obvious that $|\ell_1\cup\cdots\cup\ell_i|\le 1+(n_1-1)+\cdots+(n_i-1)$. Putting $i=t$ we get $j(L)\le 1+\mu(L)-i(L)$. If $L$ is not c-simple, i.e. $s(L)>1$, by summing up this formula evidently generalizes to $j(L)\le s(L)+\mu(L)-i(L)$.
 
 As to (b), if $L$ is c-simple and acyclic then the inequality above becomes  $|\ell_1\cup\cdots\cup\ell_i|= 1+(n_1-1)+\cdots+(n_i-1)$ for all $i$. Hence also
  $j(L)= s(L)+\mu(L)-i(L)$. If $L$ is c-simple and cyclic, then there must be a line $\ell_{i+1}$ which cuts $\ell_1\cup\cdots\cup\ell_i$ in at least {\it two} points. Consequently
  $|\ell_1\cup\cdots\cup\ell_{i+1}|< 1+(n_1-1)+\cdots+(n_{i+1}-1)$, hence $j(L)< 1+\mu(L)-i(L)$. Similarly cyclic $L$'s with $s(L)>1$ have $j(L)< s(L)+\mu(L)-i(L)$.
  
   {\bf Exercise 9C.} Show brute-force that the cyclic lattice  in Fig.9.3(A) is locally acyclic.
   
   {\bf Exercise 9D.} Show that the only BOL $\B$ of $L=\La(\Z_2\times \Z_2\times\Z_2)$ is isomorphic to the Fano plane in Fig.5.1.
   
   {\bf Exercise 9E.} Identify a triangle configuration in the BOL of Ex.9D, and accordingly (see proof of Thm.9.3) point out a cyclic localization $\B(a,1)$.
 
 {\bf Exercise 9F.} Show that the five fat lines below {\it cannot} constitute the localization $\B(a,1)$ of a BOL $\B$ whose universe of $\B(a)$ consists of the squared points.
 
 \begin{center}
 	\includegraphics[scale=0.9]{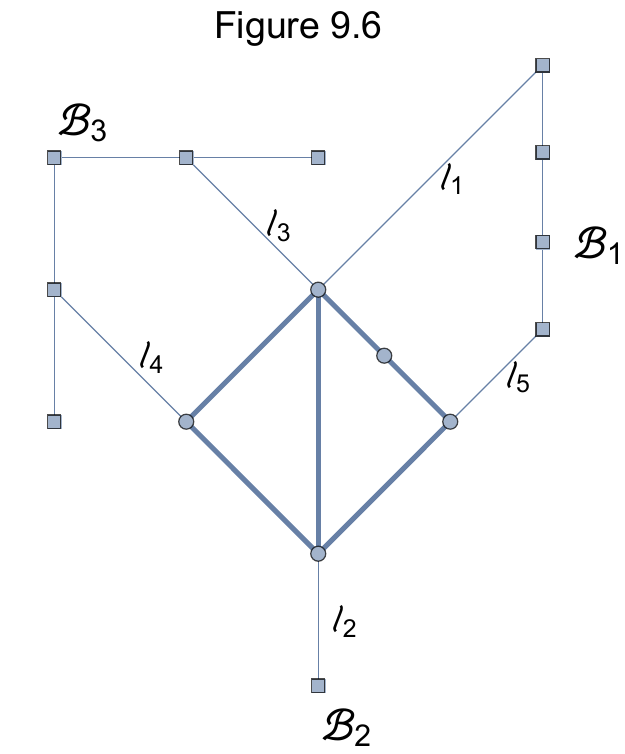} 
 \end{center}

{\it Solution to Ex.9F.}  In accordance with (2), there are lines $\ell_1',\ell_2',\ell_3'$ of $\B(a,1)$ that point to $\B_1,\B_2,\B_3$ respectively, and which form an acyclic line-set. However, $\ell_5',\ell_2',\ell_3'$ also point to $\B_1,\B_2,\B_3$, yet $T=\{\ell_1',\ell_3',\ell_5'\}$ is cyclic. Conclusion: No BOL $\B$ of a modular lattice can behave as in Fig.9.6.

 {\bf Exercise 9G.} Let $\{\ell_1,\ell_2,\ell_3\}$ be a cycle in a base of lines $\B$ of a modular lattice. (Cycles with 3 lines e.g. exist in the Fano plane). Show that the line-tops   $\ol{\ell_1},\ol{\ell_2},\ol{\ell_3}$ {\it cannot} be mutually comparable.
 
  {\it Solution to Ex.9G. }
 Let $p,q,r$ be the junctions occuring in the cycle. By assumption the $M_n$-elements $p+q\ ,p+r,\ q+r$  are mutually comparable, yet distinct since by definition of "base of lines" $\ell_i\neq \ell_j\Ra \ol{\ell_i}\neq\ol{\ell_j}$.
 We may thus assume that $p+q<p+r<q+r$. But then $q<p+q<p+r$ yields the contradiction $q+r\le p+r$.

 {\bf Exercise 9H.} Let $(x,y,z)$ be an arbitrary $M_n$-cycle of length three. Show that $x,y,z$ {\it must} be  mutually comparable. (Hence by 9G our $M_n$-cycle 'stands alone', i.e. cannot be induced by a $\B$-cycle.)
 
 {\it Solution to Ex.9H.}
Upon relabelling we can assume that $x>^*y$, say $(x_0,x_1)\searrow (y_1,y)$. 

\begin{center}
 	\includegraphics[scale=0.96]{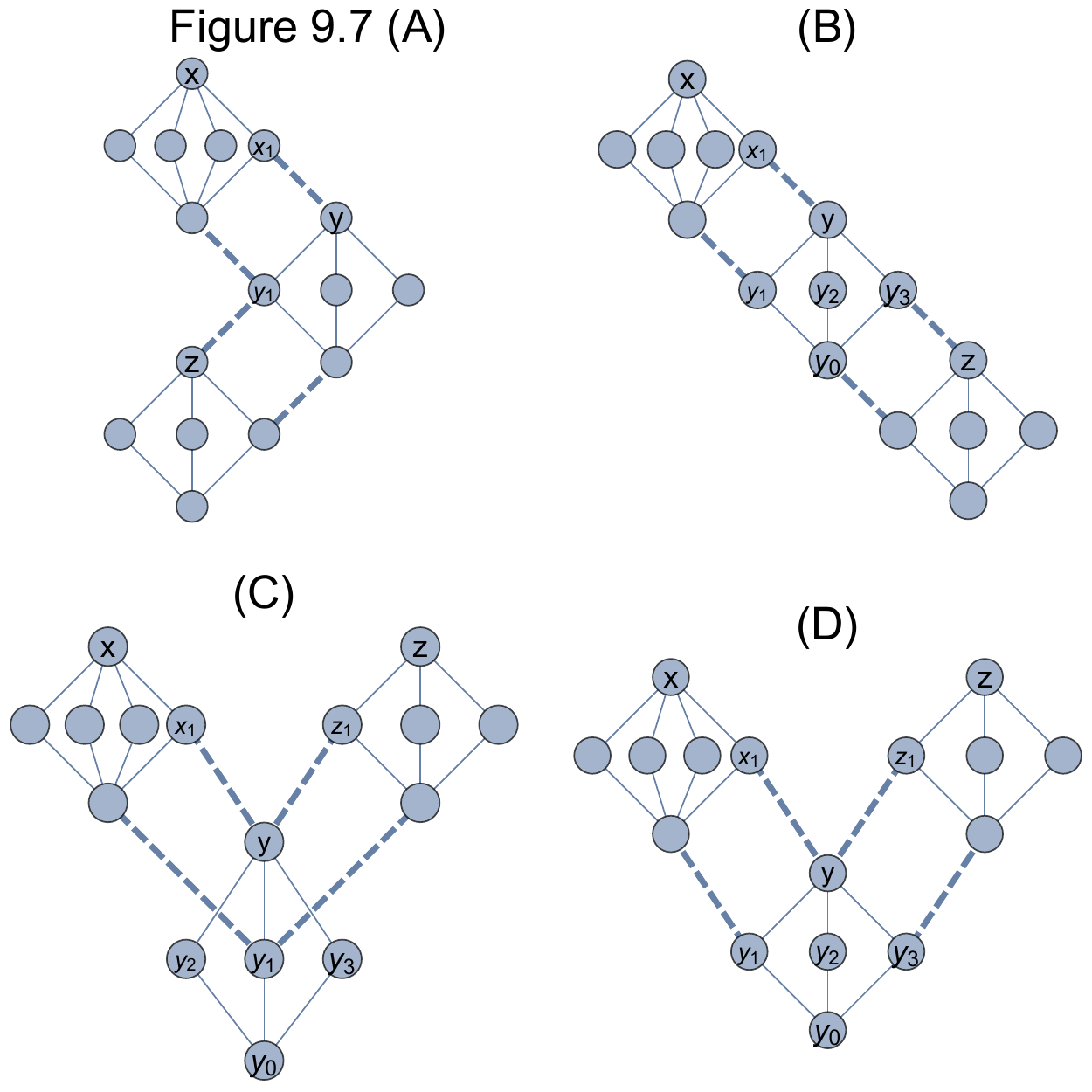} 
\end{center}

{\it Case 1:} $y>^*z$. Depending on which covering $(y_0,y_i)$ transposes down towards $z$, there are two subcases to consider. {\it Subcase 1.1:} $(y_0,y_1)\searrow (z_j,z)$ for some
$z_j\prec z$ (see Figure 9.7 (A)). This yields no $M_n$-cycle since $z<x_0$.
 {\it Subcase 1.2:} W.l.o.g. $(y_0,y_3)\searrow (z_j,z)$  (see Figure (B)), and so $x>y>z$ are mutually comparable.

 {\it Case 2:} $y<^*z$. 
 {\it Subcase 2.1:}  $(y_1,y)\nearrow (z_0,z_1)$  (see Figure (C)). Since by assumption $(x,y,z)$ is a $M_n$-cycle, we either have $y<^*z<^*x$ or $y<^*z>^*x$.
 By symmetry we can assume that $x<^*z$. The assumption that $x\le z_i$ for some $i>1$ yields the contradiction $y\le z_1x\le z_1z_i=z_0$. Hence $x\le z_1$.
   So $(z_0,z_1)\searrow (x_1,x)$, whence $x_1\le z_0$, whence  the contradiction $y\le z_0$.
 {\it Subcase 2.2:} Say $(y_3,y)\nearrow (z_0,z_1)$  (see Figure (D)).
Because $x,z$ are ccomparable, by symmetry we may assume that $x<^*z$.
If we had $x\le z_j\ (j>1)$, then $y\le z_1x\le z_1 z_j=z_0$, which is impossible. Hence $x\le z_1$. But then $y<x<z$ are mutually comparable.
 
 {\bf Exercise 9I.} Let $L=L(\Z_2^3)$. We know that the only BOL of $L$ is isomorphic (as PLS) to the Fano plane in Fig.5.1. Pick any co-atom $a\prec 1$. As stated in (2) check that indeed $|\LL(a,1)|\ge s_a$ and $j(a,1)\ge s_a+1$. Nevertheless
 $|\LL(a,1)|=6$, whereas  $j(a,1)=4$.

\vspace{2cm}
\centerline{\bf References (to be updated)}
\medskip

\begin{itemize}
	\item[[A]] M. Aigner, Combinatorial Theory, Springer Verlag Berlin, Heidelberg,
	New York 1979.
	\item[[B]] G. Birkhoff, Lattice Theory, AMS Coll. Publ., vol.25, third
	edition 1967.
	\item[[BB]] L.M. Batten, A. Beutelsbacher, The theory of finite linear spaces, Cambridge University Press 2009.
	\item[[BC]] D. Benson, J. Conway, Diagrams for modular lattices, J. Pure Appl. Algebra 37 (1985) 111-116.
	\item[[CD]] J. Crawley and R.P. Dilworth, Algebraic Theory of Lattices,
	Prentice Hall, Englewood Cliffs 1973.
	\item[[D]] A. Day, Characterizations of lattices that are bounded homomorphic
	images or sublattices of free lattices, Can. J. Math. 31 (1979) 69-78.
	\item[[DF]] A. Day and R. Freese, The role of glueing constructions in 
	modular lattice theory, p.251-260 in: The Dilworth Theorems, Birkhauser 1990.
	\item[[DP]] BA Davey, HA Priestley, Introduction to lattices and order, Cambridge University Press 2002.
	\item[[EHK]]  Marcel Ern´e, Jobst Heitzig, and J¨urgen Reinhold, On the number of distributive lattices,
	Electron. J. Combin. 9 (2002).
	\item[[EMSS]]  Erdős, P; Mullin, R. C.; Sós, V. T.; Stinson, D. R.;
	Finite linear spaces and projective planes.
	Discrete Math. 47 (1983), no. 1, 49–62.
	\item[[FH]] U. Faigle and C. Herrmann, Projective geometry on partially
	ordered sets, Trans. Amer. Math. Soc. 266 (1981) 319-332.
	\item[[FRS]] U. Faigle, G. Richter and M. Stern, Geometric exchange
	properties in lattices of finite length, Algebra Universalis 19 (1984)
	355-365.
	\item[[G]] G. Gratzer, General Lattice Theory: Foundation, Birkhauser 2011.
	\item[[H]] C. Herrmann, S-verklebte Summen von Verbaenden, Math. Zeitschrift 130 (1973) 225-274.
	\item[[HPR]] C. Herrmann, D. Pickering and M. Roddy, Geometric description
	of modular lattices, to appear in Algebra Universalis.
	\item[[HW]] C. Herrmann and M. Wild, Acyclic modular lattices and their
	representations, J. Algebra 136 (1991) 17-36.
	\item[[Hu]] A. Huhn, Schwach distributive Verbande I, Acta Sci. Math. Szeged
	33 (1972) 297-305.
	\item[[JN]] B. Jonsson and J.B. Nation, Representations of $2$-distributive
	modular lattices of finite length, Acta Sci. Math. 51 (1987) 123-128.
	\item[[K]] J. Kohonen, Generating modular lattices of up to 30 elements. Order 36 (2019) 423–435.
	\item[[KNT]] McKenzie, Ralph ; McNulty, George ; Taylor, Walter F;
	Algebras, lattices, varieties. Vol. I.
	Wadsworth and Brooks/Cole Advanced Books and Software, Monterey, CA, 1987. xvi+361 pp. 
	\item[[LMR]] 
	Lux, Klaus; Müller, Jürgen; Ringe, Michael;
	Peakword condensation and submodule lattices: an application of the MEAT-AXE. 
	J. Symbolic Comput. 17 (1994), no. 6, 529–544.
	\item[[P]] W. Poguntke,  Zerlegung von S-Verbänden, Math. Z. 142 (1975), 47–65.
	\item[[S]] M. Stern, Semimodular lattices, Teubner Texte zur Mathematik 125
	(1991), Teubner Verlagsgesellschaft Stuttgart, Leipzig, Berlin.
	\item[[Wh1]] N. White (ed.), Theory of Matroids, Encyclopedia Math. and Appl.
	26, Cambridge University press 1986.
	\item[[Wh2]] N. White (ed.), Combinatorial geometries,
	Encyclopedia Math. and Appl. 29, Cambridge University press 1987.

	\item[[W1]] M. Wild, Dreieckverbande: Lineare und quadratische 
	Darstellungstheorie, Thesis University of Zurich 1987.
	
	\item[[W2]] M. Wild, Modular lattices of finite length, first version (1992) of the present paper, ResearchGate.
	\item[[W3]] M. Wild, Cover preserving embedding of modular lattices into
	partition lattices, Discrete Mathematics 112 (1993) 207-244.
	\item[[W4]] M. Wild, The minimal number of join irreducibles of a finite modular lattice. Algebra Universalis 35 (1996),  113–123.
	
	\item[[W5]] M. Wild, Optimal implicational bases for finite modular lattices. Quaest. Math. 23 (2000), no. 2, 153–161.
	\item[[W6]] M. Wild, Output-polynomial enumeration of all fixed-cardinality ideals of a poset, respectively all fixed-cardinality subtrees of a tree, Order 31 (2014) 121–135.
	\item[[W7]] M. Wild,The joy of implications, aka pure Horn formulas: mainly a survey. Theoret. Comput. Sci. 658 (2017), part B, 264–292.
	\item[[W8]] M. Wild, Tight embedding of modular lattices into partition lattices: progress and program, Algebra Universalis 79 (2018) 1-49.
	
	\item[[Wi1]] R. Wille, Über modulare Verbände, die von einer endlichen halbgeordneten Menge frei erzeugt werden. 
	Math. Z. 131 (1973), 241–249.
	\item[[Wi2]] Subdirekte Produkte vollständiger Verbände. (German) J. Reine Angew. Math. 283(284) (1976), 53–70.
\end{itemize}

\end{document}